\documentclass[a4paper]{article}

\usepackage{amsfonts,amssymb,amsthm}

\catcode`\@=11\relax
\newwrite\@unused
\def\typeout#1{{\let\protect\string\immediate\write\@unused{#1}}}
\def\@nnil{\@nil}
\def\@empty{}
\def\@psdonoop#1\@@#2#3{}
\def\@psdo#1:=#2\do#3{\edef\@psdotmp{#2}\ifx\@psdotmp\@empty \else
    \expandafter\@psdoloop#2,\@nil,\@nil\@@#1{#3}\fi}
\def\@psdoloop#1,#2,#3\@@#4#5{\def#4{#1}\ifx #4\@nnil \else
       #5\def#4{#2}\ifx #4\@nnil \else#5\@ipsdoloop #3\@@#4{#5}\fi\fi}
\def\@ipsdoloop#1,#2\@@#3#4{\def#3{#1}\ifx #3\@nnil 
       \let\@nextwhile=\@psdonoop \else
      #4\relax\let\@nextwhile=\@ipsdoloop\fi\@nextwhile#2\@@#3{#4}}
\def\@tpsdo#1:=#2\do#3{\xdef\@psdotmp{#2}\ifx\@psdotmp\@empty \else
    \@tpsdoloop#2\@nil\@nil\@@#1{#3}\fi}
\def\@tpsdoloop#1#2\@@#3#4{\def#3{#1}\ifx #3\@nnil 
       \let\@nextwhile=\@psdonoop \else
      #4\relax\let\@nextwhile=\@tpsdoloop\fi\@nextwhile#2\@@#3{#4}}
\def\psdraft{
        \def\@psdraft{0}
}
\def\psfull{
        \def\@psdraft{100}
}
\psfull
\newif\if@prologfile
\newif\if@postlogfile
\newif\if@noisy
\def\pssilent{
        \@noisyfalse
}
\def\psnoisy{
        \@noisytrue
}
\psnoisy
\newif\if@bbllx
\newif\if@bblly
\newif\if@bburx
\newif\if@bbury
\newif\if@height
\newif\if@width
\newif\if@rheight
\newif\if@rwidth
\newif\if@clip
\newif\if@verbose
\def\@p@@sclip#1{\@cliptrue}
\def\@p@@sfile#1{
                   \def\@p@sfile{#1}
}
\def\@p@@sfigure#1{\def\@p@sfile{#1}}
\def\@p@@sbbllx#1{
                \@bbllxtrue
                \dimen100=#1
                \edef\@p@sbbllx{\number\dimen100}
}
\def\@p@@sbblly#1{
                \@bbllytrue
                \dimen100=#1
                \edef\@p@sbblly{\number\dimen100}
}
\def\@p@@sbburx#1{
                \@bburxtrue
                \dimen100=#1
                \edef\@p@sbburx{\number\dimen100}
}
\def\@p@@sbbury#1{
                \@bburytrue
                \dimen100=#1
                \edef\@p@sbbury{\number\dimen100}
}
\def\@p@@sheight#1{
                \@heighttrue
                \dimen100=#1
                \edef\@p@sheight{\number\dimen100}
}
\def\@p@@swidth#1{
                \@widthtrue
                \dimen100=#1
                \edef\@p@swidth{\number\dimen100}
}
\def\@p@@srheight#1{
                \@rheighttrue
                \dimen100=#1
                \edef\@p@srheight{\number\dimen100}
}
\def\@p@@srwidth#1{
                \@rwidthtrue
                \dimen100=#1
                \edef\@p@srwidth{\number\dimen100}
}
\def\@p@@ssilent#1{ 
                \@verbosefalse
}
\def\@p@@sprolog#1{\@prologfiletrue\def\@prologfileval{#1}}
\def\@p@@spostlog#1{\@postlogfiletrue\def\@postlogfileval{#1}}
\def\@cs@name#1{\csname #1\endcsname}
\def\@setparms#1=#2,{\@cs@name{@p@@s#1}{#2}}
\def\ps@init@parms{
                \@bbllxfalse \@bbllyfalse
                \@bburxfalse \@bburyfalse
                \@heightfalse \@widthfalse
                \@rheightfalse \@rwidthfalse
                \def\@p@sbbllx{}\def\@p@sbblly{}
                \def\@p@sbburx{}\def\@p@sbbury{}
                \def\@p@sheight{}\def\@p@swidth{}
                \def\@p@srheight{}\def\@p@srwidth{}
                \def\@p@sfile{}
                \def\@p@scost{10}
                \def\@sc{}
                \@prologfilefalse
                \@postlogfilefalse
                \@clipfalse
                \if@noisy
                        \@verbosetrue
                \else
                        \@verbosefalse
                \fi
}
\def\parse@ps@parms#1{
                \@psdo\@psfiga:=#1\do
                   {\expandafter\@setparms\@psfiga,}}
\newif\ifno@bb
\newif\ifnot@eof
\newread\ps@stream
\def\bb@missing{
        \if@verbose{
                \typeout{psfig: searching \@p@sfile \space  for bounding box}
        }\fi
        \openin\ps@stream=\@p@sfile
        \no@bbtrue
        \not@eoftrue
        \catcode`\%=12
        \loop
                \read\ps@stream to \line@in
                \global\toks200=\expandafter{\line@in}
                \ifeof\ps@stream \not@eoffalse \fi
                \@bbtest{\toks200}
                \if@bbmatch\not@eoffalse\expandafter\bb@cull\the\toks200\fi
        \ifnot@eof \repeat
        \catcode`\%=14
}       
\catcode`\%=12
\newif\if@bbmatch
\def\@bbtest#1{\expandafter\@a@\the#1
\long\def\@a@#1
\long\def\bb@cull#1 #2 #3 #4 #5 {
        \dimen100=#2 bp\edef\@p@sbbllx{\number\dimen100}
        \dimen100=#3 bp\edef\@p@sbblly{\number\dimen100}
        \dimen100=#4 bp\edef\@p@sbburx{\number\dimen100}
        \dimen100=#5 bp\edef\@p@sbbury{\number\dimen100}
        \no@bbfalse
}
\catcode`\%=14
\def\compute@bb{
                \no@bbfalse
                \if@bbllx \else \no@bbtrue \fi
                \if@bblly \else \no@bbtrue \fi
                \if@bburx \else \no@bbtrue \fi
                \if@bbury \else \no@bbtrue \fi
                \ifno@bb \bb@missing \fi
                \ifno@bb \typeout{FATAL ERROR: no bb supplied or found}
                        \no-bb-error
                \fi
                \count203=\@p@sbburx
                \count204=\@p@sbbury
                \advance\count203 by -\@p@sbbllx
                \advance\count204 by -\@p@sbblly
                \edef\@bbw{\number\count203}
                \edef\@bbh{\number\count204}
}
\def\in@hundreds#1#2#3{\count240=#2 \count241=#3
                     \count100=\count240        
                     \divide\count100 by \count241
                     \count101=\count100
                     \multiply\count101 by \count241
                     \advance\count240 by -\count101
                     \multiply\count240 by 10
                     \count101=\count240        
                     \divide\count101 by \count241
                     \count102=\count101
                     \multiply\count102 by \count241
                     \advance\count240 by -\count102
                     \multiply\count240 by 10
                     \count102=\count240        
                     \divide\count102 by \count241
                     \count200=#1\count205=0
                     \count201=\count200
                        \multiply\count201 by \count100
                        \advance\count205 by \count201
                     \count201=\count200
                        \divide\count201 by 10
                        \multiply\count201 by \count101
                        \advance\count205 by \count201
                     \count201=\count200
                        \divide\count201 by 100
                        \multiply\count201 by \count102
                        \advance\count205 by \count201
                     \edef\@result{\number\count205}
}
\def\compute@wfromh{
                \in@hundreds{\@p@sheight}{\@bbw}{\@bbh}
                \edef\@p@swidth{\@result}
}
\def\compute@hfromw{
                \in@hundreds{\@p@swidth}{\@bbh}{\@bbw}
                \edef\@p@sheight{\@result}
}
\def\compute@handw{
                \if@height 
                        \if@width
                        \else
                                \compute@wfromh
                        \fi
                \else 
                        \if@width
                                \compute@hfromw
                        \else
                                \edef\@p@sheight{\@bbh}
                                \edef\@p@swidth{\@bbw}
                        \fi
                \fi
}
\def\compute@resv{
                \if@rheight \else \edef\@p@srheight{\@p@sheight} \fi
                \if@rwidth \else \edef\@p@srwidth{\@p@swidth} \fi
}
\def\compute@sizes{
        \compute@bb
        \compute@handw
        \compute@resv
}
\def\psfig#1{\vbox {
        \ps@init@parms
        \parse@ps@parms{#1}
        \compute@sizes
        \ifnum\@p@scost<\@psdraft{
                \if@verbose{
                        \typeout{psfig: including \@p@sfile \space }
                }\fi
                \special{ps::[begin]    \@p@swidth \space \@p@sheight \space
                                \@p@sbbllx \space \@p@sbblly \space
                                \@p@sbburx \space \@p@sbbury \space
                                startTexFig \space }
                \if@clip{
                        \if@verbose{
                                \typeout{(clip)}
                        }\fi
                        \special{ps:: doclip \space }
                }\fi
                \if@prologfile
                    \special{ps: plotfile \@prologfileval \space } \fi
                \special{ps: plotfile \@p@sfile \space }
                \if@postlogfile
                    \special{ps: plotfile \@postlogfileval \space } \fi
                \special{ps::[end] endTexFig \space }
                \vbox to \@p@srheight true sp{
                        \hbox to \@p@srwidth true sp{
                                \hss
                        }
                \vss
                }
        }\else{
                \vbox to \@p@srheight true sp{
                \vss
                        \hbox to \@p@srwidth true sp{
                                \hss
                                \if@verbose{
                                        \@p@sfile
                                }\fi
                                \hss
                        }
                \vss
                }
        }\fi
}}
\catcode`\@=12\relax

\textwidth=16cm
\textheight=23cm
\oddsidemargin = -0.7truecm
\topmargin =-1truecm

\newcommand{\dr}{\partial}
\newcommand{\R}{{\bf R}}
\newcommand{\N}{{\bf N}}
\newcommand{\Z}{{\bf Z}}
\newcommand{\II}{I\hspace{-0.1cm}I}
\newcommand{\III}{I\hspace{-0.1cm}I\hspace{-0.1cm}I}
\newcommand{\tr}{\mbox{tr}}
\newcommand{\ric}{\mbox{ric}}
\newcommand{\SO}{\mbox{SO}}
\newcommand{\can}{\mbox{can}}
\newcommand{\argth}{\mbox{argth}}
\newcommand{\deltab}{\overline{\delta}}
\newcommand{\gammab}{\overline{\gamma}}
\newcommand{\eb}{\overline{e}}
\newcommand{\Gb}{\overline{G}}
\newcommand{\Kb}{\overline{K}}
\newcommand{\Pb}{\overline{P}}
\newcommand{\Sb}{\overline{S}}
\newcommand{\ricb}{\overline{\ric}}
\newcommand{\Db}{\overline{D}}
\newcommand{\Rb}{\overline{R}}
\newcommand{\omegab}{\overline{\omega}}

\newcommand{\phit}{\tilde{\phi}}
\newcommand{\Sigmat}{\tilde{\Sigma}}

\newcommand{\cM}{{\mathcal M}}

\newcommand{\bM}{{\bf M}}

\newcommand{\Rr}{\stackrel{\circ}{R}}

\newtheorem{prop}{Proposition}[section]
\newtheorem{lemma}[prop]{Lemma}
\newtheorem{conj}[prop]{Conjecture}
\newtheorem{thm}[prop]{Theorem}
\newtheorem{cor}[prop]{Corollary}
\newtheorem{remark}[prop]{Remark}

\newtheorem{df}[prop]{Definition}

\newcommand{\pg}[1]{\paragraph{#1}}

\newenvironment{thn}[1]{\vskip 0.2cm \noindent{\bf Theorem #1.} \it}{\rm
\vspace{0.2cm}} 
\newenvironment{lmn}[1]{\vskip 0.2cm \noindent{\bf Lemma #1.} \it}{\rm
\vspace{0.2cm}} 

\newcommand{\btm}{\begin{thm}}
\newcommand{\etm}{\end{thm}}
\newcommand{\blm}{\begin{lemma}}
\newcommand{\elm}{\end{lemma}}
\newcommand{\bcr}{\begin{cor}}
\newcommand{\ecr}{\end{cor}}
\newcommand{\bdf}{\begin{df}}
\newcommand{\edf}{\end{df}}
\newcommand{\bprop}{\begin{prop}}
\newcommand{\eprop}{\end{prop}}
\newcommand{\bas}{\begin{asser}}
\newcommand{\eas}{\end{asser}}
\newcommand{\beq}{\begin{equation}}
\newcommand{\eeq}{\end{equation}}
\newcommand{\bpv}{\begin{proof}}
\newcommand{\epv}{\end{proof}}
\newcommand{\bit}{\begin{itemize}}
\newcommand{\eit}{\end{itemize}}
\newcommand{\bpn}{\begin{pfn}}
\newcommand{\epn}{\end{pfn}}
\newcommand{\btn}{\begin{thn}}
\newcommand{\etn}{\end{thn}}
\newcommand{\bln}{\begin{lmn}}
\newcommand{\eln}{\end{lmn}}

\newenvironment{pfn}[1]{\vskip 0.2cm \noindent{\it Proof #1.}}{$\square$
\vspace{0.2cm}}

\newcommand{\Sib}{\overline{\Sigma}}
\newcommand{\gb}{\overline{g}}
\newcommand{\cC}{\mathcal{C}}
\newcommand{\cD}{\mathcal{D}}
\newcommand{\Met}{\mathcal{M}et}
\newcommand{\Imm}{\mathcal{I}mm}
\newcommand{\CMet}{\mathcal{CM}et}
\newcommand{\CImm}{\mathcal{CI}mm}
\newcommand{\gab}{\overline{\gamma}}
\newcommand{\hyp}{\mathbf{H}^3}
\newcommand{\dhyp}{\partial\hyp}
\newcommand{\cL}{\mathcal{L}}
\newcommand{\isom}{\mathrm{Isom}}
\newcommand{\db}{\overline{\partial}}

\begin{document}

\title{Hypersurfaces in $H^n$ and the space of
its horospheres}

\author{Jean-Marc Schlenker\thanks{
Laboratoire Emile Picard, UMR CNRS 5580,
Universit{\'e} Paul Sabatier,
118 route de Narbonne,
31062 Toulouse Cedex 5,
France.
\texttt{schlenker@picard.ups-tlse.fr; http://picard.ups-tlse.fr/\~{
}schlenker}. }}

\date{January 2001}

\maketitle

\begin{abstract}

A classical theorem, mainly due to Aleksandrov \cite{Al} and
Pogorelov \cite{Po}, states that any Riemannian metric on $S^2$ with
curvature $K>-1$ is induced on a unique convex surface in $H^3$. A
similar result holds with the induced metric replaced by the third
fundamental form. We show that the same phenomenon happens with yet
another metric on immersed surfaces, which we call the horospherical
metric. 

This results extends in higher dimension, the metrics obtained are then
conformally flat. One can also study equivariant immersions of surfaces
or the metrics obtained on the boundaries of hyperbolic
3-manifolds. Some statements which are difficult or only conjectured for
the induced metric or the third fundamental form become fairly easy when
one considers the horospherical metric. 

The results concerning the third fundamental form are obtained using a
duality between $H^3$ and the de Sitter space $S^3_1$. In the same way,
the results concerning the horospherical metric are proved through a
duality between $H^n$ and the space of its horospheres, which is
naturally endowed with a fairly rich geometrical structure.

\bigskip

\begin{center} {\bf R{\'e}sum{\'e}} \end{center}

Un th{\'e}or{\`e}me bien connu, d{\^u} essentiellement {\`a} Aleksandrov \cite{Al} et
Pogorelov \cite{Po}, affirme que chaque m{\'e}trique {\`a} courbure $K>-1$ sur
$S^2$ est induite sur une unique surface convexe dans $H^3$; un r{\'e}sultat
analogue est valable lorsque la m{\'e}trique induite est remplac{\'e}e par la
troisi{\`e}me forme fondamentale. On montre ici que le m{\^e}me ph{\'e}nom{\`e}ne se
produit si on consid{\`e}re une autre m{\'e}trique sur les surface, qu'on
appelle m{\'e}trique horosph{\'e}rique. 

Ce r{\'e}sultat s'{\'e}tend en dimension plus grande, les m{\'e}triques obtenues
{\'e}tant alors conform{\'e}ment plates. On peut aussi {\'e}tudier les immersions
{\'e}quivariantes de surfaces ou les m{\'e}triques
obtenues sur les bords de vari{\'e}t{\'e}s hyperboliques de dimension $3$, et
des {\'e}nonc{\'e}s difficiles ou seulement conjectur{\'e}s pour la m{\'e}trique
induite ou la troisi{\`e}me forme fondamentale deviennent faciles pour cette
m{\'e}trique horosph{\'e}rique. 

Les r{\'e}sultats concernant la troisi{\`e}me forme fondamentale sont obtenus en
utilisant une dualit{\'e} connue entre $H^3$ et l'espace de Sitter $S^3_1$;
de la m{\^e}me mani{\`e}re, les r{\'e}sultats concernant la m{\'e}trique horosph{\'e}rique
d{\'e}coulent d'une dualit{\'e} qu'on d{\'e}crit entre $H^n$ et l'espace de ses
horosph{\`e}res, qui est muni naturellement d'une structure g{\'e}om{\'e}trique
assez riche. 

\end{abstract}

\pg{Convex surfaces in $H^3$}
Let $S$ be a smooth, strictly convex, compact surface in $H^3$. Then $S$ is
diffeomorphic to $S^2$, and the Gauss formula indicates that its induced
metric has curvature $K>-1$. A well-known theorem, to which several
mathematicians have contributed (e.g. Weyl, Nirenberg \cite{N},
Aleksandrov \cite{Al,AZ} and Pogorelov \cite{Po}; see \cite{L1} for a
modern approach):

\btm \label{thm-h3}
Each smooth metric with curvature $K>-1$ on $S^2$ is induced on a unique
convex surface in $H^3$. 
\etm

Note that a similar result holds in $\R^3$, and also in the
3-dimensional sphere $S^3$. The uniqueness here is of course up to
global isometries of $H^3$. 

Although the ``usual'' way of considering this theorem is as a statement
on surfaces in $H^3$, it can also be understood as a remarkable
statement of existence and uniqueness for a strongly non-linear boundary
value problem: finding a hyperbolic metric on the 3-dimensional ball
$B^3$ which induces a given metric on the boundary. When considered in
this way a basic question is whether the boundary condition chosen here
is the only one possible, or indeed the best. One of the goals of this
paper is to show that there is an alternative candidate.

\pg{The third fundamental form of a surface}
This is a fairly classical bilinear form, called $\III$ here, on the
tangent of an immersed  surface ($\III$ is defined in section 1). 
When the surface is strictly convex, it
is a Riemannian metric; for surfaces in $\R^3$, $\III$ is just the
pull-back by the Gauss map of the canonical metric on $S^2$. 

An interesting point is that $\III$ provides another good boundary
condition for the existence and uniqueness of hyperbolic metrics on
$B^3$:

\btn{\cite{these}} \label{thm-a3}
Let $h$ be a smooth metric on $S^2$. $h$ is the third fundamental form
of a convex surface $S$ in $H^3$ if and only if it has curvature
$K<1$. $S$ is then unique up to global isometries. 
\etn

This result is quite strongly related to analoguous polyhedral
statements -- just like theorem \ref{thm-h3} was related to the
investigation of polyhedra in $H^3$ (see \cite{alex}). See \cite{hr,dap}
for some related questions.

\pg{The de Sitter space}
Theorem \ref{thm-a3} is strongly related to an interesting duality
between $H^3$ and the 3-dimensional de Sitter space, which will be
denoted here $S^3_1$. $S^3_1$ is the 3-dimensional, geodesically
complete, simply connected Lorentz space with constant curvature $1$. It
is described more fully in section 1. The duality between $H^3$ and
$S^3_1$ -- also briefly recalled in section 1 -- associates to each
point of $H^3$ a totally geodesic, space-like plane in $S^3_1$, and to
each totally geodesic (oriented) plane in $H^3$ a point in $S^3_1$. Thus
each strictly convex surface $S$ in $H^3$ has a well defined ``dual surface''
$S^d$ in $S^3_1$, which is space-like and convex. Moreover, the third
fundamental of $S$ is the induced metric on $S^d$. Theorem \ref{thm-h3}
can therefore be considered as an isometric embedding theorem in $S^3_1$
-- and indeed that is how it is proved.

\pg{The space of horospheres in $H^n$}
Just like $S^3_1$ is the space of oriented planes in $H^3$ with a
natural geometric structure, we can consider the space of horospheres in
$H^n$. This is done in some details in section 1; this space -- called
$C^n_+$ here for reasons that should become clear -- is $n$-dimensional,
and it has a natural
degenerate metric, of signature $(n-1,0)$. It also has a natural
foliation by curves, which we call ``vertical lines'', which are tangent
to the kernel of the degenerate metric, and come with a canonical
parametrization. One of the goals of this paper is to show that one can
do some interesting geometry in this space which, in a sense, can be
considered as a ``degenerate'' constant curvature space. 

Moreover, there is a natural duality between $H^n$ and $C^n_+$, sending
a hypersurface in $H^n$ to the set $S^*$ of points in $C^n_+$ corresponding to
the horospheres tangent to $S$; this duality -- on which some details
are given in section 2 -- has some similarities with the $H^3/S^3_1$
duality.

\pg{The ``horospherical metric'' of a surface}
When $S$ is a smooth oriented hypersurface in $H^n$, we define the
``horospherical metric'' on $S$ as the ``metric'' induced on the dual
``hypersurface'' $S^*$, and we denote it by $I^*$. In general $I^*$
might be degenerate or otherwise badly behaved, just like the third
fundamental form of a non-convex surface in $H^3$. There is a natural
class of hypersurfaces in $H^n$, however, for which $I^*$ is a smooth
Riemannian metric: the ``H-convex'' hypersurfaces, which are those which
at each point lie on one side of their tangent horosphere. If $S$ is
such an H-convex hypersurface, its dual $S^*$ is a smooth,
``space-like'' hypersurface in $C^n_+$, and moreover it is convex in a
natural sense (see section 3). 

The metric $I^*$ has a simple expression in terms of the usual
extrinsic invariants of a hypersurface: $I^*=I+2\II+\III$, where $\II$
is the second fundamental form. The point is that it provides another
good boundary condition for the existence and uniqueness of hyperbolic
metrics on $B^3$. There is class of metrics on $S^2$, which we call
``H-admissible'' (resp. ``C-admissible''), and which have a rather simple
definition (see 
definition \ref{df-adm}); those metrics have curvature $K<1$
(resp. $K\in (-1, 1)$), and are exactly the horospherical metrics of the
H-convex (resp. convex) surfaces in $H^3$:  

\btn{\ref{ii3}}
Let $h$ be a smooth metric on $S^2$. It is the horospherical metric
$I^*$ of a H-convex immersed sphere $S$ in $H^3$ if and only if it is
H-admissible. It is the horospherical metric of a convex embedded sphere
$S\subset H^3$ if and only if it is C-admissible. In each case, $S$ is
unique up to the global isometries of $H^3$.  
\etn

In higher dimension, it is not so clear what the metrics induced on
e.g. the convex hypersurfaces are. Of course not all metrics are
possible, and the conformal flatness of the metrics plays a role
\cite{cartan}. It turns out that the situation is much simpler for the
horospherical metric, since here again a simple results holds and is
easy to prove. 

\btn{\ref{iin}}
Let $h$ be smooth metric on $S^{n-1}$. $h$ is the horospherical metric
$I^*$ of a H-convex sphere $S$ in $H^n$ if and only if: 
\begin{itemize}
\item $h$ is conformal to $\can_{S^{n-1}}$;
\item $h$ is H-admissible, in the sense that it is conformal to
$\can_{S^{n-1}}$ and that 
$ 2\ric_h - \frac{S_h}{n-2} - (n-3) h $ is everywhere negative definite. 
\end{itemize}
$S$ is then unique up to the isometries of $H^n$.
Moreover, $S$ is convex if and only if all eigenvalues of
$2(n-2)\ric_h-S_h h$ are in $(-(n-2)(n-3), (n-2)(n-3))$. 
\etn

Again, there is also a simple characterization of the metrics which are
the horospherical metrics of convex hypersurfaces. 

\pg{Equivariant surfaces}
Let $\Sigma$ be a surface of genus at least two. Although $\Sigma$
carries many metrics with curvature $K>-1$, they can of course not be
induced by an embedding in $H^3$, since it should then be convex. One needs the
slightly refined notion of equivariant embedding. That is a couple
$(\phi, \rho)$, where $\phi$ is an embedding of the universal cover
$\Sigmat$ of $\Sigma$, and $\rho$ is a morphism from $\pi_1(\Sigma)$
into $\isom(H^3)$, such that:
$$ \forall x\in \Sigmat, \forall \gamma\in \pi_1(\Sigma), \phi(\gamma
x)=\rho(\gamma)\phi(x)~. $$

One can then search for equivariant embeddings inducing a given metric;
it turns out that (because of the index theorem) there are too many of
those, so that one can impose an additionnal condition on $\rho$. 

\begin{thm}[Gromov \cite{PDR}]
Let $\Sigma$ be a surface of genus at least $2$, and let $h$ be a smooth
metric on $\Sigma$ with curvature $K>-1$. There is an equivariant
isometric embedding $(\phi, \rho)$ of $(\Sigma, h)$ into $H^3$ such that
$\rho$ fixes a plane. 
\end{thm}
A remarkable point is that it is still not know whether the uniqueness
holds in the theorem above. On the other hand, an analoguous results
holds with the induced metric replaced by the third fundamental form:

\begin{thm}[\cite{iie}] \label{tiie}
Let $\Sigma$ be a surface of genus at least $2$, and let $h$ be a smooth
metric on $\Sigma$ with curvature $K<1$. There is a unique equivariant
embedding $(\phi, \rho)$ of $\Sigma$ into $H^3$ such that the third
fundamental form $\III$ of $\phi$ is $h$ and that $\rho$ fixes plane. 
\end{thm}

The uniqueness above is of course up to global isometries of $H^3$. 

Considering the horospherical metric on $\Sigma$ instead of either the
induced metric or the third fundamental form leads to simpler results
again. There are simple definitions (see \ref{df-hcvx-g}) of
``H-admissible'' and ``C-admissible'' metrics on $\Sigma$, which are
sub-classes of the  metrics with curvature $K<1$ and $K\in (-1,1)$
respectively. Then:

\btn{\ref{ie}}
A smooth metric $h$ on $\Sigma$ is the horospherical metric of a
H-convex equivariant immersion whose representation fixes a plane if and
only if $h$ is H-admissible. It is the horospherical metric of a convex
embedding whose representation fixes a plane if and
only if $h$ is C-admissible. The equivariant immersion/embedding is then
unique up to global isometries. 
\etn

Here again results also hold in higher dimension, and the proof is quite
simple.

\pg{Manifolds with boundaries}
As stated above, theorem \ref{thm-h3} can be considered as a boundary
value problem for hyperbolic metrics on the 3-dimensional ball. When
considered in 
this way it should be possible to generalize it to manifolds other than
$B^3$. Such a generalization was proposed in the following conjecture. 

\begin{conj}[Thurston] \label{thurston}
Let $M$ be a 3-dimensional manifold with boundary which admits a
complete, convex co-compact metric. Then, for any smooth metric $h$ on $\dr
M$ with curvature $K>-1$, there is a unique hyperbolic metric $g$ on $M$
which induces $h$ on the boundary, and for which the boundary is convex. 
\end{conj}

The proof of the existence part of the conjecture was
obtained by Labourie \cite{L5,L4}, but the uniqueness remains
unknown. Theorem \ref{tiie} also suggests that the same kind of result
might hold with the induced metric replaced by the third fundamental
form; actually the main point of this paper is that the "horospherical
metric" works quite well for this kind of results. There are natural
classes of ``H-admissible'' and ``C-admissible'' metrics, defined in
\ref{h-adm-bord}, which have curvature $K<1$, such that:

\btn{\ref{thm-bord}}
Let $h$ be a smooth metric on $\dr M$. 
\begin{enumerate}
\item $h$ is the horospherical metric
of a H-convex immersion $\phi$ of $\dr M$ in $M$ for a complete
hyperbolic metric  $g$ on $M$, such that the image of $\phi$ can be deformed
through immersions to the boundary at infinity of $M$, if and only if
$h$ is H-admissible. $g$ and $\phi$ are then unique.  
\item $h$ is the horospherical metric of $\dr M$ for a hyperbolic metric
$g$ on $M$, such that $\dr M$ is convex, if and only if $h$ is
C-admissible. $g$ is then unique. 
\end{enumerate}
\etn

In this setting again, the proof is easy, although it uses a deep
result, the Ahlfors-Bers theorem (seen here as a bijection between
conformal structures on $\dr M$ and hyperbolic metrics on $M$; see
\cite{ahlfors}). Actually, theorem \ref{ie} is a direct consequence of
theorem \ref{thm-bord}; it might still be helpful to some readers to
have stated it separately.

\pg{Hypersurfaces in $S^n_1$}
Note that if $S\subset H^n$ is a convex surface, and if $S^d$ is the
dual surface in $S^n_1$ (which is also a convex surface, and moreover is
space-like) then the first, second and third fundamental forms on $S^d$
are $I^d=\III$, $\II^d=\II$, and $\III^d=I$ respectively. Therefore:
$$ I + 2\II + \III = I^d + 2\II^d + \III^d~, $$
so that most of the themes described in this paper for
hypersurfaces in $H^n$ are also valid for convex hypersurfaces in
$S^n_1$, and can be proved by considering the dual surface in
$H^n$. Presumably a weaker hypothesis than convexity could be used (like
the H-convexity condition in $H^n$); it should be possible to repeat
some of the arguments below without reference to the dual hypersurface
in $H^n$, by replacing the horospheres in $H^n$ by their dual
hypersurfaces in $S^n_1$.

\paragraph{} 
The main point of all this is that some
results which are either rather difficult or actually still conjectures
for the induced metric or the third fundamental form of (hyper-)surfaces
become easy when one considers the horospherical metric instead. Section
7 contains examples of some other areas where this metric
might be of interest.

\section{The space of horospheres in $H^n$}

We will describe in this section the natural geometric structure on the
space of horocycles in $H^n$. This structure will be a basic tool in the
sequel, so it will be important to understand various basic aspects of
it, for instance what the ``hyperplanes'' or the ``umbilical
hypersurfaces'' are, and how the isometries act.

\pg{Horospheres in $H^n$}
To any point at infinity $y\in \dr_\infty H^n$, one associates a
``Buseman function'', defined, up to the addition of a constant, as: 
$$ B_y(x) = \lim_{t\rightarrow \infty} d(\gamma(t), x_0) - d(\gamma(t),
x)~, $$
where $x_0\in H^n$ is any fixed point and $\gamma$ is any geodesic ray
with $\lim_{t\rightarrow \infty}\gamma(t)=y$. The level sets of the
Buseman functions are called the horospheres of $H^n$. By construction,
each horosphere is associated to a unique point at infinity; if two
horospheres have the same point at infinity, then they are equidistant. 

Let $H$ be a smooth oriented hypersurface in $H^n$, and let $X$ and $Y$ be two
vector fields on $H$. Call $D$ the Levi-Civit{\`a} connection of
$H^n$. Then:
$$ D_XY = \Db_XY + \II(X,Y)N~, $$
where $\Db$ is the Levi-Civit{\`a} connection of the induced metric $I$ on
$H$, and $N$ is the unit normal vector field on $H$. $\II$ is called the
{\bf second fundamental form} of $H$, it is a symmetric bilinear form on
$TH$. The {\bf Weingarten operator} $B$ of $H$ is then defined by:
$$ \II(X,Y) = I(-BX,Y) = I(X,-BY)~. $$
The sign convention used here is not so standard but will make things
easier because we will want to use the exterior normal of e.g. spheres
in $H^n$. 
The {\bf third fundamental form} of $H$ is:
$$ \III(X,Y) = I(BX,BY)~. $$
Horospheres in $H^n$ are characterized by the equation: 
$$ \III = \II = I~, $$
in particular they are umbilical.

\pg{The Poincar{\'e} model of $H^n$}
There is a convenient model of hyperbolic $n$-space, called the Poincar{\'e}
model, which is a conformal map from $H^n$ to the Euclidean disc
$D^n$ (see e.g. \cite{GHL}). Moreover, there is a also a conformal map
from the $n$-sphere 
$S^n$ minus a point to the Euclidean space $\R^n$. This map can be
obtained by stereographic projection. Composing those maps gives us a
conformal map from $H^n$ to a geodesic ball in $S^n$, whose radius can
be chosen by choosing the right radius for the image of the Poincar{\'e}
model of $H^n$ in $\R^n$.

\pg{The Klein model of $H^{n+1}$ and the $H^{n+1}$/$S^{n+1}_1$ duality}
There is also another model of $H^{n+1}$, called the ``Klein'' or
``projective'' model. This is a map from $H^{n+1}$ to $D^{n+1}$ which
has the striking property that the geodesics of $H^{n+1}$ are mapped to
the segments of $D^{n+1}$. 

It has a natural extension to a projective model
of a ``hemisphere'' of the $n+1$-dimensional de Sitter space $S^{n+1}_1$
on the complement of $D^{n+1}$ in $\R^{n+1}$. 
$S^{n+1}_1$ can be also seen as a
quadric in Minkowski $n+2$-space, with the induced metric, as follows: 
$$ S^{n+1}_1 = \{ x\in \R^{n+2}_1 ~ | ~ \langle x,x\rangle = 1 \}~. $$
There is a natural duality
between $H^{n+1}$ and $S^{n+1}_1$, which associates to a totally
geodesic hyperplane in $H^{n+1}$ a point in $S^{n+1}_1$. It can be
defined in the Minkowski models of $H^{n+1}$ and $S^{n+1}_1$ as
follows. Remember that $H^{n+1}$ can be seen as:
$$ H^{n+1} = \{ x\in \R^{n+2}_1 ~ | ~ \langle x,x\rangle = -1  ~
\mbox{and} ~ x_0 >0 \}~. $$
Given a point $x\in H^{n+1}$, let $D$ be the line going through $0$ and
$x$ in $\R^{n+2}_1$, and let $D^d$ be its orthogonal, which is a
space-like hyperplane in $\R^{n+2}_1$. The dual of $x$ is the
intersection $x^d := D^d\cap S^{n+1}_1$. The same works in the opposite
direction, from points in $S^{n+1}_1$ to oriented hyperplanes in
$H^{n+1}$.

Given a smooth, oriented, strictly convex hypersurface $S\subset H^{n+1}$, 
the set of points in $S^{n+1}_1$
which are the duals of the hyperplanes tangent to $S$ is called the dual
hypersurface; the notations used here will be $S^d$. It is a space-like,
convex hypersurface in $S^{n+1}_1$; 
its induced metric is $=I^d=\III$, while its
third fundamental form is $\III^d=I$. 
See e.g. \cite{shu,rh} for a detailed construction and some additional
remarks (in particular concerning polyhedra) on the projective model.

\pg{The space of horospheres in $H^n$}
The Poincar{\'e} model of $H^n$ therefore allows us to consider $H^n$ as the
interior of a geodesic sphere $S_0$ in $S^n$. $S^n$ can, through the Klein
model of $H^{n+1}$,
be considered as the boundary at infinity of $H^{n+1}$. 
$S_0$ is the boundary of a totally geodesic hyperplane $H_0\subset
H^{n+1}$; let $S_0^*$ be the point in $S^{n+1}_1$ which is dual of $H_0$. 
The horospheres
in $H^n$ are then identified with the spheres in $S^n$ which are
interior to and tangent to $S_0$; they are the boundaries of the
totally geodesic hyperplanes in $H^{n+1}$ which have a point at infinity
in $S_0$. The set of point in $S^{n+1}_1$ which are the duals of those
hyperplanes is part of the cone of lines in $\R^{n+1}$ going through
$S_0^*$ and tangent to $S^n$; more precisely, it is the set of points of
this cone which lie strictly between $S_0^*$ and $S^n$, or, in other
terms, the positive light-cone of a point in $S^{n+1}_1$ -- whence the
notation $C^n_+$. 

We already see that $C^n_+$ inherits from this construction a degenerate
metric -- the one induced on the cone by the de Sitter metric -- and a
foliation by a family of lines -- those going through $S_0^*$. We call
those lines ``{\bf vertical}''. By
construction both the metric and the family of vertical lines are
independant of the choices made in the construction. The vertical lines
are actually characterized as the curves which are everywhere tangent to
the kernel of the (degenerate) metric $g_0$. 

Note that $C^n_+$ has, by construction, a very large group of
"isometries" which fix both $g_0$ and the vertical lines. This indicates
that it is a kind of "degenerate constant curvature space". 

\vspace{0.3cm}
\centerline{\psfig{figure=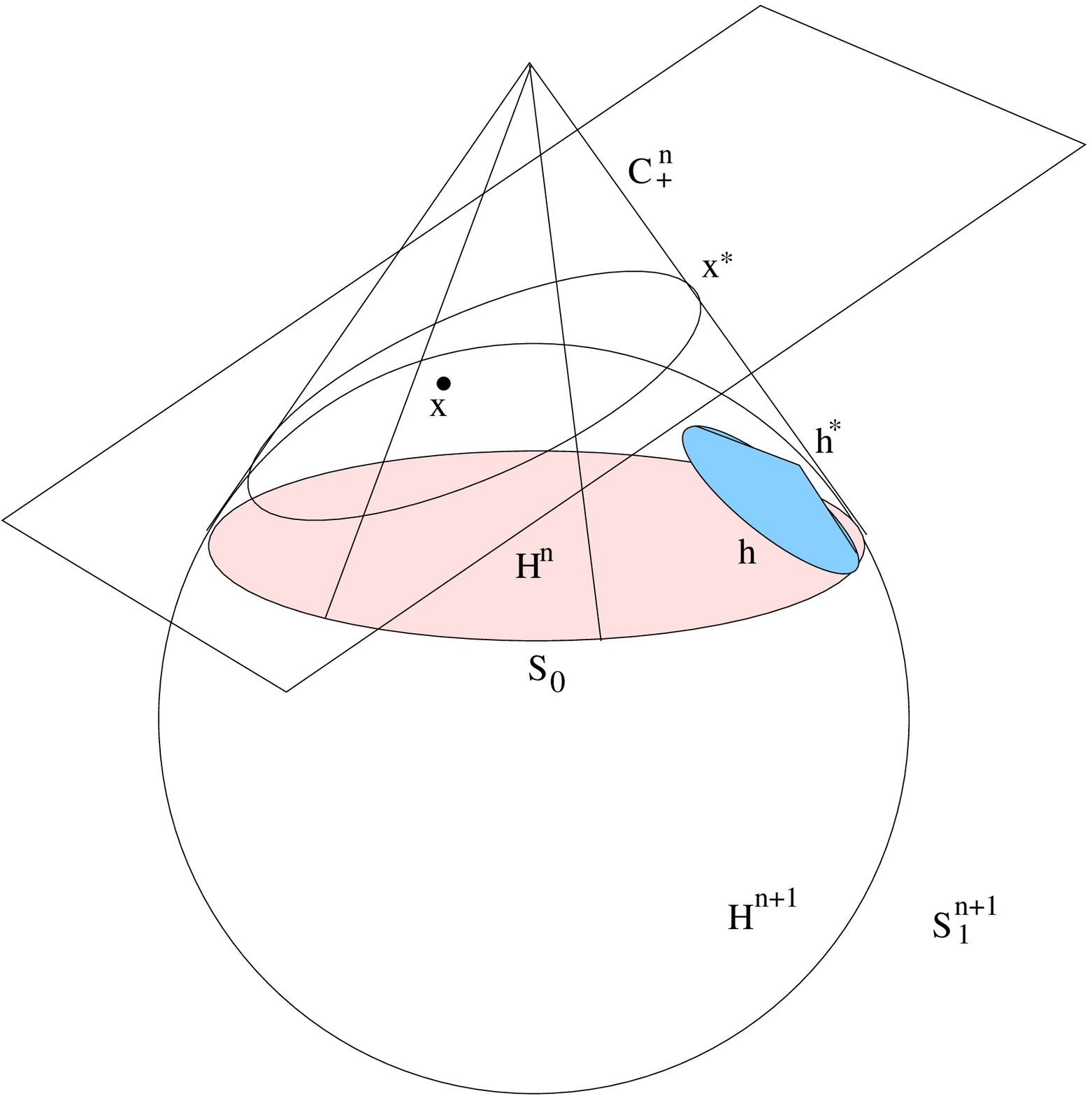,height=8cm}} \vspace{0.3cm} 
\centerline{\bf Figure 1: conical model of $C^n_+$} \vspace{0.3cm}

\pg{A cylindrical model}
A slightly different model, which might sometimes be more convenient, is
obtained by taking $H^n$ as a hemisphere in $S^n$; $S_0$ is then an
``equatorial'' $n-1$-sphere, and its dual point $S_0^*$ is at infinity,
so that $C^n_+$ is identified with the union of the lines tangent to
$S^n$ at a point of $S_0$, and parallel to the line in $R^{n+1}$ which
is orthogonal to the hyperplane containing $S_0$. 

\vspace{0.3cm}
\centerline{\psfig{figure=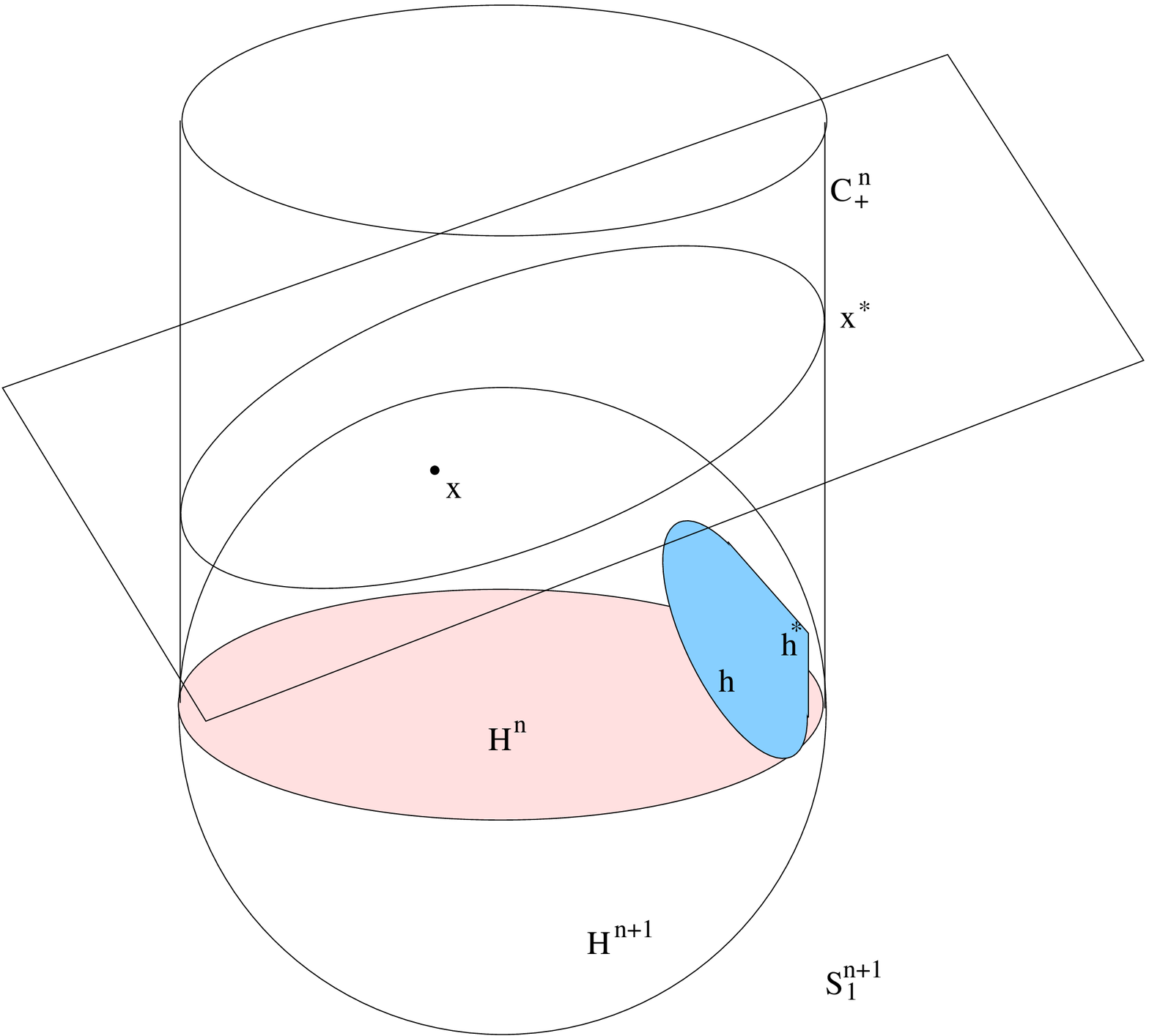,height=8cm}} \vspace{0.3cm} 
\centerline{\bf Figure 2: cylindrical model of $C^n_+$} \vspace{0.3cm}

\pg{The induced structure}
As a submanifold of $S^{n+1}_1$, $C^n_+$ inherits a degenerate
metric $g_0$, i.e. a bilinear form on the tangent space which is at each
point of rank $n-1$. Moreover the kernel of this bilinear form, which at
each point is made of a line in the tangent space, integrate as
``lines'' in $C^n_+$. 

Those lines are the lines in $\R^{n+1}$ which contain $S_0^*$ and are
tangent to $S^n$. They are therefore light-like geodesics of
$S^{n+1}_1$, and are naturally equiped with a connection; in other terms
they have a parametrization by $\R$ which is defined up to an affine
transformation. But those lines actually also have a natural
parametrization (up to a constant); namely, it is easy to check that
they correspond to the sets of horospheres which have a given focal
point at infinity, so that the horospheres
corresponding to two points in a given line are equidistant. The distance
between them defines the required parametrization. 

Note that, from $g_0$ and this canonical parametrization of the vertical
lines, one could define a (family of) Riemannian metrics on $C^n_+$. But
it does not seem very helpful to do this.

\pg{Totally geodesic hyperplanes}
$C^n_+$ comes equiped with a collection of hypersurfaces which play a
special role, and that will be called ``{\bf totally geodesic
hyperplanes}''. They are the sets of points duals to the horospheres
containing a given point in $H^n$. In the cone model described above, they
correspond to the intersections of the cone with the hyperplanes of
$\R^{n+1}$ which are tangent to $S^n$ at an interior point of $S_0$. 
Thus the metric induced on those totally geodesic hyperplanes is
isometric to the canonical metric on $S^{n-1}$. Moreover, it should be
clear from the description below that they are the only space-like
hypersurfaces in $C^n_+$ with an induced metric isometric to $(S^{n-1},
\can)$. 

By definition, the set of those totally geodesic hyperplanes is an
$n$-dimensional manifold -- it is parametrized by $H^n$. 

\blm \label{tangent}
Let $x\in C^n_+$, and let $P\subset T_xC^n_+$ be a 
hyperplane which is transverse to the vertical line at $x$. 
There is a unique totally geodesic hyperplane $H_0$ in $C^n_+$ which is
tangent to $P$ at $x$. 
\elm

\bpv 
Consider the cylindrical model of $C^n_+$ described above. $P$
corresponds to an $n-1$-plane in $\R^{n+1}$ which is disjoint from
$S^n$. There are two hyperplanes containing $P$ which are tangent to
$S^n$, and one of them is tangent to $C^n_+$ along a line; so there is a
unique hyperplane $\Pb$ which contains $P$, is 
transverse to $C^n_+$, and is tangent to $S^n$. $\Pb$ intersects $C^n_+$
along an $n-1$-dimensional manifold which, by construction, is a totally
geodesic hyperplane in $C^n_+$. 
\epv

\pg{Parallel transport along the vertical lines}
In the cone model above, the tangent space to $C^n_+$ is parallel (in
$S^{n+1}_1$) along the ``vertical lines'' (which are the lines in
$C^n_+$ which are tangent to $S^n$ at the points of $S_0$). Therefore, the
restriction of the Levi-Civit{\`a} connection of $S^{n+1}_1$ defines a
connection along the vertical lines in $C^n_+$, and thus also a natural
notion of parallel transport along those lines. We call this induced
connection $D^V$.

\pg{A kind of connection}
Now let $x_0\in C^n_+$, and let $H$ be a hyperplane in $T_xC^n_+$ which
is transverse to the vertical direction. We can define a kind
of connection, which we call $D^H$, along the vectors tangent to $H$ at
$x_0$. Note that it depends on the choice of $H$ ! 
It is defined as follows. Call $H_0$ the totally geodesic
hyperplane tangent to $H$ at $x_0$, let $X\in H$, and let $Y$ be a
vector field defined in a neighborhood of $x_0$, which is tangent to
$H_0$; then define:  
$$ D^H_X Y = D^0_X Y~, $$
where $D^0$ is the Levi-Civit{\`a} connection of $H_0$ for the induced
metric. Moreover, if $T$ is the vector field everywhere parallel to the
vertical lines, and with length given by the natural parametrization of
those lines, then we decide that, for any function $f$ on $C^n_+$:
$$ D^H_X fT = df(X) T~. $$
This clearly defines $D^H_XY$ by linearity for any vector field $Y$ on
$C^n_+$.  

Note that, on the other hand, we do not define a canonical connection on
$C^n_+$ -- and we will not really need one here. 

The definition of $D^H$ can also be obtained in an extrinsic way as
follows. For $x_0$ and $H$ chosen as above, there is a unique hyperplane
of $\R^{n+1}$ which is transverse to $C^n_+$, tangent to $S^n$, and
contains $H$. This plane contains a unique light-like line $D'$ containing
$x_0$. Now choose $X\in H$, and let $Y$ be a vector field defined on
$C^n_+$ in a neighborhood of $x_0$. One can project on $T_{x_0}C^n_+$
along $D'$ the vector $D_X^{S^{n+1}_1}Y$, where $D^{S^{n+1}_1}$ is the
Levi-Civit{\`a} connection of $S^{n+1}_1$. The reader might want to check
that this indeeds defines the same vector as $D^H_XY$. Of course the
point is that the result depends on $D'$, and therefore on $H$.

\section{H-convex hypersurfaces in $H^n$}

This section contains some elementary remarks about the dual, in $C^n_+$
of some hypersurfaces in $H^n$. They are then used to give an intrinsic,
and quite simple, expression of the metric on $C^n_+$.

\pg{H-convex hypersurfaces}
The following notion of convexity is important in our context. 

\bdf
Let $S$ be an oriented hypersurface in $H^n$, let $x\in S$, and let $h$
be a horosphere in $H^n$. We say that $h$ is {\bf tangent} to $S$ at $x$
if $h$ is tangent to $S$ at $x$ in the usual sense, and moreover the
convex side of $h$ is on the exterior side of $S$. 
\edf

\bdf
Let $S$ be an oriented hypersurface in $H^n$. $S$ is {\bf H-convex} if, at
each point $x\in S$, $S$ remains on the concave side of the horosphere
tangent to $S$ at $x$. $S$ is {\bf strictly H-convex} if, moreover, the
distance 
between $S$ and that horosphere does not vanish up to the second order
in any direction at $x$.  
\edf

From now on, ``H-convex'' will be understood as ``strictly H-convex''
except when otherwise stated. ``Convex'' will also mean ``strictly
convex''. 

\bdf
Let $S$ be a hypersurface in $H^n$. We denote by $S^*$ the set
of points in $C^n_+$ which are dual to the horospheres tangent to $S$. 
\edf

\bdf 
Let $S$ be a smooth hypersurface in $C^n_+$. We say that $S$ is {\bf
space-like} if $S$ is everywhere transverse to the vertical lines. 
\edf

It is not difficult to check that the only compact space-like
hypersurfaces in $C^n_+$ are spheres. 

We will often implicitely identify a hypersurface $S$ with its dual,
using the natural map sending a point $x\in S$ to the dual $h^*$ of the
horosphere tangent to $S$ at $x$. 
 
\blm \label{dual-metric}
If $S$ is an H-convex hypersurface in $H^n$ such that its principal
curvatures are nowhere equal to $-1$, then $S^*$ is an immersed
surface in $C^n_+$. This happens in particular when $S$ is H-convex, and
$S^*$ is then space-like. The metric induced by $g_0$ on $S^*$ is: 
$$ I^* := I + 2 \II + \III~, $$
where $\II$ and $\III$ are the second and third fundamental forms of $S$
respectively. 
\elm

Note that, for instance here, the identification of $I^*$ with $I + 2
\II + \III$ 
implicitely uses the map from $S$ to $S^*$ sending a point of $S$ to the
point of $S^*$ which is dual to the horosphere tangent to $S$ at that
point. 

The proof of this lemma will use the cylindrical model of $C^n_+$ in an
explicit way. Consider an H-convex hypersurface $H$ in $H^n$, and let
$x\in H$. We will use the cylindrical model of $C^n_+$, with $x$ located
at the "north pole" of
$S^n$; the dual of the horosphere $h$ which is tangent to $H$ at
$x$ is then a point $h^*$  of the intersection of $C^n_+$ (seen as a
cylinder) with the hyperplane in $\R^{n+1}$ which is tangent to $S^n$ at
$x$. 

The tangent space to
$H$ at $x$ is identified with an affine
$n-1$-dimensional subspace $V$ of $\R^{n+1}$, and the tangent space to
$C^n_+$ at $h^*$ can be seen as an $n$-dimensional affine subspace $W$ of
$\R^{n+1}$ which contain an $n-1$-plane parallel to $V$. We call $\phi$
the duality map from $H$ to $H^*$, sending a point $y$ in $H$ to the dual of
the horosphere tangent to $H$ at $y$, and we consider $d\phi$ as a
map from $V$ to $W$, where $W\supset V$. Then:

\bprop \label{prop}
The linearized map at $x$ is $T_x\phi = E+B$, where $E$ is the identity
map on $V=T_xH$.  
\eprop

\bpv
Let $v\in T_xH$; call $v^*$ the vector in $W$ corresponding to the
variation of the dual point to the horosphere tangent to $H$ at a point
which moves in the direction of $V$ on $H$. $v^*$ is the sum of a term
$v_1^*$ corresponding to the displacement of $x$ (with a parallel
transport of the tangent hyperplane) and a term $v_2^*$ corresponding to
the variation of the tangent hyperplane, while $x$ doesn{\'{}}t move. Using
the cylindrical model, one checks that $v_1^*=v$ (with both terms seen
as in $W'$) while $v_2^*=Bv$. 
\epv

\bpn{of lemma \ref{dual-metric}}
The previous proposition shows that $S^*$ is smooth except maybe
when $B$ has $-1$ as one of its eigenvalues. 

Moreover, the bilinear form induced on $W$ by $g_0$ (i.e. by the de
Sitter metric on the outside of the ball) is a degenerate metric which
coincides, on the parallel transport of $V$, with the metric induced on
$V$ by $H^n$. Therefore, if $v, v'\in T_xH$, we have that $v^*, v'^*\in
V$ and:
$$ \langle v+Bv, v'+Bv'\rangle = \langle v,v'\rangle + \langle Bv,
v'\rangle + \langle v, Bv'\rangle + \langle Bv, Bv'\rangle~, $$
so that: 
$$ \langle v^*, v'^* \rangle = I(v,v') + 2\II(v, v') + \III(v,v')~, $$
and the result follows.   
\epn

\pg{This is a duality}
An important point is that the map sending a hypersurface $S$ in $H^n$ to
its dual $S^*$ in $C^n_+$ is a real duality, in the following
sense. First remark that to each totally geodesic hyperplane $H_0$ in
$C^n_+$ is associated a point in $H^n$, namely the intersection of all
the horospheres duals to the point of $H_0$. We call this point the dual
of $H_0$, and denote it by $H_0^*$. Then:

\blm \label{duality}
If $S^*$ is smooth, then $S$ is the set of points in $H^n$ which are
duals of a totally geodesic hyperplane in $C^n_+$. 
\elm

\bpv
This follows again from proposition \ref{prop}, and from the
correspondance between vectors on $S$ and on $S^*$.
\epv

\pg{An intrinsic definition of the metric $g_0$}
The previous lemma can be used to give a simple form of the metric on
$C^n_+$; using it will relieve us from the constant use of the cone
model, the de Sitter space and so on. 

\blm \label{intrinsic}
There exists an isometry $\Phi$ from $C^n_+$ to $S^{n-1}\times \R$ with
the (degenerate) metric: 
$$ g_0 \simeq e^{2t} \can_{S^{n-1}}~, $$
where $\can_{S^{n-1}}$ is the canonical metric on $S^{n-1}$. Moreover
the vertical lines are sent to the lines $\{ s\}\times \R$, for $s\in
S^{n-1}$, with the same parametrization. 
\elm

\bpv
Let $x_0\in H^n$. For $t\in \R\setminus \{ 0\}$, call $S_t$ the geodesic
sphere of radius $|t|$ centered at $x_0$, with the normal oriented
towards the exterior for $t>0$ and towards the interior for
$t<0$. Define a map $\Psi$ from $S^{n-1}\times (\R\setminus \{ 0\})$ to
$C^n_+$ sending $(s, t)$ to the horosphere tangent to $S_t$ at the point
$\exp_{x_0} (ts)$, where $s$ is considered as a unit vector in
$T_{x_0} H^n$.  $\Psi$ can then be extended by continuity to a map from
$S^{n-1}\times \R$ to $C^n_+$. $\Phi$ is the inverse of $\Psi$. 

By lemma \ref{dual-metric}, the metric induced on $S_t^*$ is:
\begin{eqnarray} 
I^*_t & = & I_t + 2\II_t + \III_t \nonumber \\
& = & \sinh^2(t) \can_{S^{n-1}} (1 + 2 \coth(t) + \coth^2(t)) \nonumber \\
& = & (\sinh^2(t) + 2\sinh(t)\cosh(t) + \cosh^2(t)) \can_{S^{n-1}}
\nonumber \\ 
& = & e^{2t} \can_{S^{n-1}}~. \nonumber 
\end{eqnarray}

Now, using e.g. the cylindrical model described above, with $x_0$ as the
``north pole'' in $S^n$, shows that the
surfaces $S_t^*$ are the intersections of $C^n_+$ (seen as a cylinder in
$R^{n+1}$) with the horizontal hyperplanes, i.e. the hyperplanes in
$\R^{n+1}$ which are parallel to the hyperplane containing
$S_0$. Therefore the lines $\{ s\} \times \R$ are in the kernel of $g_0$,
and moreover they correspond to the vertical lines. Finally, by
definition of their parametrization (by the distance between equidistant
horospheres) it is the same as the one they have in $S^{n-1}\times 
\R$. 
\epv

\pg{A conformal map}
Now we remark that all the space-like hypersurfaces in $C^n_+$ can
be naturally identified in a conformal way; they are moreover all
naturally conformal to the boundary at infinity of $H^n$. 
Call $\Pi_0$ the map from $C^n_+$ to $\dr_\infty H^n$ sending a
horosphere to its point at infinity. Then:

\begin{lemma} \label{lm-conf}
\begin{enumerate}
\item let $H_1$ and $H_2$ be two compact space-like hypersurfaces in
$C^n_+$. The projection from $H_2$ to $H_1$ along the vertical lines is
conformal for the induced metrics on $H_1$ and $H_2$. 
\item For each space-like hypersurface $H_1\subset C^n_+$, the
restriction of $\Pi_0$ to $H_1$ is conformal for the induced metric on
$H_1$ and the usual conformal metric on $\dr_\infty H^n$. 
\end{enumerate}
\end{lemma}

\bpv
The first point is a direct consequence of lemma \ref{intrinsic}
above. 
For the second point remark that, if $x_0$ is the point in $H^n$ which
is the intersection of the horospheres in $H_1$, then the map sending a
horosphere $h\in H_1$ to its point at infinity is by construction an
isometry between $H_1$ with its induced metric and $\dr_\infty H^n$ with
the visual metric at $x_0$. It is therefore a conformal map. 
\epv

Let $H$ be an oriented hypersurface in $H^n$; there is a natural map
from $H$ to $\dr_\infty H^n$, which sends a point $x\in H$ to the end
point of the ray starting at $x$ in the direction of the oriented normal
vector to $H$ at $x$. We call this ``Gauss map'' $G$ (see e.g. \cite{L5}
for some applications of this map). 
As a consequence of lemma \ref{lm-conf} we obtain the following: 

\blm \label{conf}
If $H$ is an H-convex hypersurface in $H^n$, the conformal structure
obtained on $H$ as the pull-back by $G$ of the conformal structure on
$\dr_\infty H^n$ is the conformal structure of $I^*$. 
\elm

\pg{Umbilical hyperplanes}
Some hypersurfaces in $C^n_+$ play a special role and have a very simple
geometry; they are the surfaces $H^*$, where $H$ is an umbilical
hypersurface in $H^n$. By lemma \ref{dual-metric}, $H^*$ is then
homothetic to $H$. This is specially interesting when $H$ is a totally
geodesic hyperplane in $H^n$, since then $H^*$ is isometric to
$H$. We call those surfaces ``{\bf dual hyperplanes}''. It is not difficult to
check that the image of a dual hyperplane by the projection on a totally
geodesic hyperplane along the vertical lines is a hemisphere.

\pg{Isometries}
Let $\gamma$ be an isometry of $H^n$. Consider the cone model of $C^n_+$
described in section 1. Then $\gamma$ acts on $S^n$ as a M{\"o}bius
transformation leaving $S_0$ stable. Therefore it acts as an isometry on
$H^{n+1}$, 
seen as the interior of $S^n$, and therefore also as an isometry on the
de Sitter space which lies on the outside of $S^n$, leaving
invariant the cone made of the (light-like) lines tangent to $S^n$ along
$S_0$ and containing $S_0^*$. So, by construction, $\gamma$ also acts 
on $C^n_+$ without changing its metric or its vertical lines.

Note that if $\gamma$
has no fixed point in $\dr_\infty H^n$, then it has no fixed
point in $C^n_+$ -- since an isometry fixing a horosphere should fix its
point at infinity. This strongly contrasts with the $H^n/S^n_1$
duality, where all isometries of $H^n$ without fixed point in $H^n\cup
\dr_\infty H^n$ have at least one fixed point in $S^n_1$. 

The isometries of $C^n_+$ can be characterized in the
following simple ways. 

\blm \label{isometries}
\begin{enumerate}
\item Let $H$ be a totally geodesic hyperplane in $C^n_+$. For any
isometry $\gamma$ of $H^n$, (the extension to $C^n_+$ of) $\gamma$,
composed with the projection on $H$ 
  along the vertical lines is a conformal transformation of $H$.
\item Moreover, any conformal transformation of $H$ corresponds in this
  way to a unique isometry. 
\item Let $D$ be a dual hyperplane in $C^n_+$. Each isometry of $D$
  extends in exactly two ways as an isometry of $C^n_+$, one of which
  preserves orientation. 
\end{enumerate}
\elm

\bpv
Isometries correspond by defintion to isometries of $H^n$, which act
conformally on $\dr_\infty H^n$, and thus on $H$ by lemma
\ref{lm-conf}; point (1) follows. Conversely, any conformal
transformation of $H$ defines by lemma \ref{lm-conf} a conformal
transformation of $\dr_\infty H^n$, and therefore an isometry of $H^n$,
and also an isometry of $C^n_+$. This proves point (2). 

For point (3), let $D^*$ be the dual hyperplane of $D$, i.e. the totally
geodesic hyperplane in $H^n$ such that $D$ corresponds to the set of
horospheres tangent to $D^*$. Let $\gamma$ be an isometry of $D$. By
construction, $D^*$ is isometric to $D$, so that $\gamma$ defines an
isometry $\gamma^*$ of $D^*$. Since $D^*$ is an hyperplane in $H^n$,
$\gamma^*$ has two extensions as an isometry of $H^n$, one of which
preserves orientation. We call this orientation preserving extension
$\gammab^*$ again. 
$\gammab^*$ defines a unique isometry $\gammab$ of $C^n_+$, which
leaves $D$ stable by construction. The same works for the other
extension of $\gamma^*$. 
\epv

\pg{Group actions}
It might be interesting to understand what the quotient of
e.g. $C^2_+$ (resp. $C^3_+$) by the action of discrete group acting
co-compactly on $H^2$ (resp. $H^3$) is .

\section{More on the dual of a hypersurface}

\pg{Second fundamental forms in $C^n_+$}
Let $H$ be a hypersurface in $C^n_+$. Let
$x\in H$, and call $H_0\in T_xC^n_+$ the totally geodesic hyperplane
tangent to $H$ at $x$. Let $X$ and $Y$ be vector fields on $H$. Locally
(in the neighborhood of $x$) $H$ intersects exactly once each vertical line;
therefore, the "vertical connection" $D^v$ defined in section 1 allows
us to extend $X$ and $Y$ as vector fields on a neighborhood of $x$ in $C^n_+$
by parallel transport along the vertical lines. We can then unse the
kind of connection defined in section 1 to define a ``second fundamental
form'' of $H$ at $x$.

\bdf \label{df-II}
The second fundamental form of $H$ at $x$ is defined as:
$$ \II^*(X,Y) := \Pi (D^{H_0}_XY)~, $$
for the extended vector fiels, 
where $\Pi$ is the projection on the vertical direction in $T_xC^n_+$
along the direction of $H_0$. 
\edf

\blm \label{ptes-II}
\begin{enumerate}
\item  $\II^*$ defines a symmetric bilinear form on $H_0$. 
\item If $P_0$ is the
(unique) totally geodesic hyperplane in $C^n_+$ which is tangent to
$H_0$ at $x$, then $H$ is locally the graph of a function $u$ above
$P_0$; $\II^*$ is then the hessian of $u$ at $x$ for the metric induced
on $P_0$. 
\item $\II^*$ is also the hessian at $x$ of $u$, seen as a function on
$H$, for the induced metric $I^*$ on $H$. 
\end{enumerate}
\elm

In the second part of this lemma, $u$ is the function such that, at a
point $y\in P_0$ near $x$, $u(y)$ is the "oriented distance" from $y$ to the
intersection of $H$ with the vertical line through $y$, for the natural
parametrization of that vertical line. 

\bpv
The first point is obviously a consequence of the other. For the second
point note that, in the neighborhood of $x$, the extended vector field
$Y$ is of the form:
$$ Y = Y_0 + du(Y_0) T~, $$
with $Y_0$ tangent to $P_0$. Therefore the definition of $D^{H_0}$
shows that:
$$ D^{H_0}_XY = D^0_XY_0 + (X.du(Y_0)) T~, $$
where $D^0$ is the Levi-Civit{\`a} connection of the induced metric $g_0$ on
$P^0$, and the result follows since $du=0$ at $x$. 

For the third point note that, by lemma \ref{intrinsic},
$I^*=e^{2u}g_0$, so that the Levi-Civit{\`a} connection $D^*$ of $I^*$ is
given by: 
$$ D^*_XY = D^0_XY + du(X)Y + du(Y)X - g_0(X,Y) D^0u~, $$
where vector fields on $H$ and $P_0$ are identified through the
projection along the vertical lines. Therefore (by the usual conformal
transformation formulas, see e.g. \cite{Be}, chapter 1):
$$ (D^*du)(X,Y) = (D^0du)(X,Y) - 2du(X)du(Y) + g_0(X,Y) \|
du\|^2_{g_0}~, $$
so that $D^*du = D^0du$ at $x$ since $du=0$ at $x$. 
\epv

We then use $\II^*$ to define the "Weingarten operator" of a
hypersurface $H$ in $C^n_+$:

\bdf
If $H$ is a space-like hypersurface in $C^n_+$ and $x\in H$, the
"Weingarten operator" of $H$ at $x$ is the linear map $B^*$ from $T_xH$ to
$T_xH$, self-adjoint for $I^*$, defined by:
$$ \II^*(X,Y) = I^*(B^*X,Y) = I^*(X,B^*Y)~. $$ 
\edf

\pg{An inversion formula}
We have already seen in lemma \ref{dual-metric} that:
$$ I^*(X,Y) = I((E+B)X, (E+B)Y)~. $$
Together with the previous lemma, it shows that:

\blm \label{inversion}
If $S$ is a hypersurface in $H^n$ with no principal curvature equal to
$-1$ at any point, then:
$$ B^* = (E+B)^{-1}~. $$
\elm

\pg{Convex hypersurfaces}
Using the previous definition, we can define a convex hypersurface in
$C^n_+$:

\bdf
Let $H$ be a space-like hypersurface in $C^n_+$. We say that $H$ is {\bf
convex} if $B^*$ is positive definite at each point
of $H$. $H$ is {\bf tamely convex} if all eigenvalues of $B^*$ are in
$(0,1)$ at each point. 
\edf

The point is that convex hypersurfaces in $C^n_+$ have smooth dual
hypersurfaces in $H^n$, and that tamely convex hypersurfaces have convex
duals. More precisely:

\blm \label{trans-inverse}
Let $H$ be a hypersurface in $C^n_+$ such that $B^*$ is nowhere
degenerate. Then $H^*$ is smooth, and its induced metric is:
$$ I(X,Y) = I^*(B^*X, B^*Y)~.  $$
$H$ is tamely convex if and only if $H^*$ is convex. 
\elm

\bpv
This follows again from proposition \ref{prop}. 
\epv

\section{Isometric embeddings in $C^n_+$}

The point of this section is to give an elementary study of the induced
metrics on hypersurfaces in $C^n_+$, like the one which can be found in
elementary differential geometry books for hypersurfaces in
e.g. $\R^n$. The results are a little different, however, due to the
degeneracy of the metric.

\pg{The Gauss formula}
The curvature tensor of the induced metric on a hypersurface in $C^n_+$
is determined by the following analog of the Gauss formula:

\blm \label{gauss}
Let $H$ be a space-like hypersurface in $C^n_+$. Let $x\in H$, call
$P_0$ the (unique) totally geodesic hyperplane in $C^n_+$ which is
tangent to $H$ at $x$. Let $X,Y,Z$ be three vector fields on $H$. The
Riemann curvature tensor $R^*$ of the induced metric $I^*$ on $H$ is
given by:  
$$ R^*_{X,Y} Z = R^0_{X,Y}Z + \II^*(X,Z)Y -\II^*(Y,Z) X - I^*(Y,Z)B^*X
+ I^*(X,Z) B^*Y~, $$
where $R^0$ is the curvature tensor of $P_0$.
\elm

Note that this formula differs from the Euclidean one, in particular
because it is linear in $B^*$ instead of quadratic. 

\bpv
We call also $X,Y$ and $Z$ the projections of the vector fields on
$P_0$, and $g_{P_0}$ its metric, which has constant curvature $1$. The
metric on $H$ is then the pull-back of $e^{2u}g_{P_0}$ under the
projection of $H$ to $P_0$ along the vertical lines. Therefore, the
Levi-Civit{\`a} connection $\Db$ of $I^*$ is (see e.g. \cite{Be}, chap. 1): 
$$ \Db_XY = D_XY + du(X)Y + du(Y)X - g_{P_0}(X,Y) Du~, $$
where $D$ is the Levi-Civit{\`a} connection of $g_{P_0}$. Thus, using the
fact that $du=0$ at $x$, we find that, still at $x$: 
\begin{eqnarray}
R^*_{X,Y}Z & = & \Db_X\Db_Y Z - \Db_Y\Db_XZ - \Db_{[X,Y]}Z \nonumber \\
& = & D_X\Db_Y Z - D_Y\Db_XZ - D_{[X,Y]}Z \nonumber \\
& = & D_X(D_Y Z + du(Y)Z + du(Z)Y - g_{P_0}(Y,Z)Du) - \nonumber \\
&& - D_Y(D_X Z + du(X)Z + du(Z)X - g_{P_0}(X,Z)Du) - D_{[X,Y]}Z \nonumber \\
& = & R^0_{X,Y}Z + (D_Xdu)(Y)Z + (D_Xdu)(Z)Y - \nonumber \\ 
&& - (D_Ydu)(X)Z - (D_Ydu)(Z)X
- I^*(Y,Z) D_XDu + I^*(X,Z)D_YDu~, \nonumber
\end{eqnarray}
and the result follows.
\epv

\pg{Some consequences}
To simplify somewhat the exposition, we concentrate here on surfaces,
i.e. the $n=3$ case. The above formula becomes, for the Gauss curvature
of a surface:
$$ K^* = 1 - \tr(B^*)~. $$
From lemma \ref{inversion}, this can be translated as:
$$ K^* = 1 - \tr((E+B)^{-1}) = 1 - \frac{\tr(E+B)}{\det(E+B)}~, $$
so that:
$$ K^* = \frac{\det(E+B) - \tr(E+B)}{\det(E+B)} = \frac{\det(B)-1}{1 +
\tr(B) + \det(B)}~, $$
and, since $\det(B)-1$ is the Gauss curvature $K$ of the dual surface by
the (usual) Gauss formula in $H^3$:
$$ K^* = \frac{K}{K+2H+2}~, $$
where $H$ is the mean curvature of the dual surface in $H^3$. 
Therefore, when $K\neq 0$, we have: 
$$ K^* = \frac{1}{1+2(H+1)/K}~. $$
Thus the constant mean curvature $-1$ surface in $H^3$ are characterized
as those whose dual has constant curvature $1$ (of course the minus sign
is just a question of orientation).

\pg{The Codazzi theorem}
Another basic point is that, just as for hypersurfaces in Euclidean
space, we have:

\blm \label{codazzi}
Let $H$ be a space-like hypersurface in $C^n_+$ with a smooth dual
hypersurface, and let $D^*$ be the
Levi-Civit{\`a} connection of its induced metric. Then, for any vector fields
$X, Y$ on $H$:
$$ D^*_XY = B^* D_X(B^{* -1}Y)~, $$
and:
$$ (D^*_X B^*)Y = (D^*_YB^*)X~. $$
\elm

\bpv
For the first part of the lemma, we want to show that the connection
(again called $D^*$) defined by:
$$ D^*_XY = (E+B)^{-1}D_X((E+B)Y) $$
is torsion-free and compatible with $I^*$. But it is torsion-free
because:
\begin{eqnarray}
D^*_XY - D^*_YX & = & (E+B)^{-1} (D_X ((E+B)Y)-D_Y(E+B)X) \nonumber \\
& = & (E+B)^{-1} ((E+B)(D_XY-D_YX) + \nonumber \\
&& + (D_XE)Y - (D_YE)X + (D_XB)Y -
(D_YB)X) \nonumber \\ 
& = & D_XY-D_YX~, \nonumber
\end{eqnarray}
the last step using the Codazzi equation on the dual
hypersurface. Therefore $D^*_XY - D^*_YX = [X,Y]$, and $D^*$ is
torsion-free. 

To check that $D^*$ is compatible with $I^*$ is even more direct. If
$X,Y,Z$ are vector fields on $H$, then:
\begin{eqnarray}
X.I^*(Y,Z) & = & X.I((E+B)Y, (E+B)Z) \nonumber \\
& = & I(D_X((E+B)Y), (E+B)Z) + I((E+B)Y, D_X((E+B)Z)) \nonumber \\
& = & I^*(D^*_XY,Z) + I^*(Y,D^*_XZ)~.  \nonumber 
\end{eqnarray}

The second point of the lemma is easy to prove using the first; if $X$
and $Y$ are vector fields on $H$, then:
\begin{eqnarray}
(D^*_X B^*)Y - (D^*_YB^*)X & = & D^*_X (B^*Y) - D^*_Y (B^*X) -
B^*(D^*_XY - D^*_YX) \nonumber \\
& = & B^* (D_XY-D_YX) - B^*(D^*_XY - D^*_YX) \nonumber \\
& = & B^*[X,Y] - B^*[X,Y] \nonumber \\
& = & 0~. \nonumber 
\end{eqnarray}
\epv

\pg{Remark}
Lemma \ref{codazzi} provides another proof of the formulas given above,
relating $K$ and $K^*$ for surfaces in $H^3$ and in $C^3_+$. Indeed,
let $(e_1, e_2)$ be an orthonormal frame on a surface $S\subset H^3$;
then, by definition of $I^*$, $(\eb_1,
\eb_2):=((E+B)^{-1}e_1,(E+B)^{-1}e_2)$ is an orthonormal frame for $I^*$
on $S^*$. Moreover, the connection 1-forms $\omega$ and $\omegab$ of
those frames are the same:
\begin{eqnarray}
\omega(u) & := & I(D_ue_1, e_2) \nonumber \\
& = & I^*((E+B)^-1D_ue_1, (E+B)^-1e_2) \nonumber \\
& = & I^*(D^*_u ((E+B)^-1e_1), (E+B)^-1e_2) \nonumber \\
& = & I^*(D^*_u \eb_1, \eb_2) \nonumber \\
& =: & \omegab(u)~. \nonumber 
\end{eqnarray}
Therefore, the curvatures on $S$ and $S^*$ differ only by the same
factor as the area forms, so that:
$$ K^* = \frac{K}{\det(E+B)}~. $$

\pg{Induced metrics -- higher dimensions}
Here we take $n\geq 4$, the next paragraph will center on $n=3$. Let $h$
be a smooth metric on $S^{n-1}$, we have the following elementary
characterization of whether $h$ can be obtained as the induced metric on
a space-like hypersurface in $C^n_+$. 

\btm \label{isom-n}
$(S^{n-1}, h)$ admits a space-like isometric embedding into $C^n_+$ if
and only if $h$ is conformal to $\can_{S^{n-1}}$. In this case the
embedding is unique up to the isometries of $C^n_+$.
\etm

\bpv
Let $P_0$ be any totally geodesic hyperplane in $C^n_+$. If $(S^{n-1},
h)$ has a space-like isometric embedding in $C^n_+$, then the projection
from the image to $P_0$ is conformal by lemma \ref{lm-conf}. Therefore
$h$ is conformal to $\can_{S^{n-1}}$. Conversely, if $h$ is conformal to
$\can_{S^{n-1}}$ then there exists a function $u:S^{n-1}\rightarrow \R$
such that $h=e^{2u}\can_{S^{n-1}}$; then the graph of $u$ above $P_0$
is, by lemma \ref{intrinsic}, isometric to $h$. 
\epv

A more interesting -- but still easy -- question is to determine when
$h$ is induced on a convex or tamely convex hypersurface in $C^n_+$. We
call $S_h$ the scalar curvature of $h$. 

\btm \label{cvx-n}
$h$ is induced on a convex space-like hypersurface $H$ in $C^n_+$ if and
only if $h$ is conformal to $\can_{S^{n-1}}$ and 
$ 2\ric_h - \frac{S_h}{n-2} - (n-3) h $ is everywhere negative definite.
$H$ is then unique up to isometries of $C^n_+$. $H$ is tamely convex
if and only if all eigenvalues of $2(n-2)\ric_h-S_h h$ are in
$(-(n-2)(n-3), (n-2)(n-3))$. 
\etm

We will say that $h$ is {\bf H-admissible} if it satisfies the
"positive definite" hypothesis of the theorem, and that $h$ is {\bf
C-admissible} if it satisfies the eigenvalue hypothesis. 

\bpv
Let $(e_i)_{1\leq i\leq n-1}$ be an orthonormal frame for $I^*$ which
diagonalizes $B^*$, and let $(k_i)_{1\leq i\leq n-1}$ be the associated
eigenvalues of $B^*$. Call $K_{i,j}$ the sectional curvature of $h$ on
the 2-plane generated by $e_i$ and $e_j$. Then, by lemma \ref{gauss}:
$$ K_{i,j} = 1 - k_i - k_j~, $$
so that:
$$ \ric_h(e_i,e_i) = \sum_{j\neq i} K_{i,j} 
= (n-2) -(n-3)k_i - \sum_j k_j~,  $$
and
$$ S_h = \sum_i \ric_h(e_i,e_i) =  (n-1)(n-2) - 2(n-2) \sum_j k_j~, $$
so that:
\begin{eqnarray}
k_i & = & \frac{S_h + (n-2)(n-3) - 2(n-2)\ric_h(e_i,e_i)}{2(n-2)(n-3)}
\nonumber \\ 
& = & \frac{S_h - 2(n-2)\ric_h(e_i,e_i)}{2(n-2)(n-3)} + \frac{1}{2}~,
\nonumber 
\end{eqnarray}
and both results follow. 
\epv

\pg{Induced metrics -- $n=3$}
The analog of theorem \ref{isom-n} is even simpler in dimension $n=3$,
since in that case all metrics on $S^2$ are conformal to the canonical
metric. Therefore, if $h$ is a smooth metric on $S^2$: 

\btm \label{isom-3}
$(S^2, h)$ admits a unique (up to the isometries of $C^3_+$) isometric
embedding in $C^3_+$.
\etm

To understand the metrics induced on convex surfaces we have to
introduce a definition (which is also a lemma).

\bdf \label{df-adm}
Let $h$ be a smooth metric on $S^2$. Let $x\in S^2$. There is a unique
function $u_x$ on $S^2$ such that the metric $e^{-2u_x}h$ has constant
curvature $1$ and that $u_x(x)=du_x(x)=0$. We say that $h$ is
{\bf H-admissible} if, for each $x\in S^2$, the hessian of $u_x$ at $x$ is
positive definite, and that $h$ is {\bf C-admissible} if, for each $x$,
all eigenvalues of the hessian of $u_x$ at $x$ are in $(0,1)$. 
\edf

\bpv
We have to prove the existence and uniqueness of $u_x$. 

$h$ is conformal to $\can_{S^2}$, so there exists a function
$u:S^2\rightarrow \R$ such that $e^{2u}\can_{S^2} = h$. Choose a totally
geodesic plane $P_0\subset C^3_+$, and let $S$ be the graph of $u$ above
$P_0$. Then, by lemma \ref{intrinsic}, the metric induced on $S$ is
$h$. 

Now let $x\in S$. By lemma \ref{tangent}, there exists a unique totally
geodesic plane $P_1$ in $C^3_+$ which is tangent to $S$ at $x$. $P_1$ is
the graph above $S$ of a function $v$ on $S$. Then $e^{-2v}h$ is the
metric induced on $P_1$, and is isometric to $\can_{S^2}$, so $v$
satisfies the conditions set on $u_x$. 

Conversely, if $w:S\rightarrow \R$ satisfies those conditions, then the
graph $P$ of $w$ above $S$ has as induced metric $\can_{S^2}$, so it is
a totally geodesic plane, and moreover it is tangent to $S$ at
$x$. Thus, by lemma \ref{tangent}, $P=P_1$, and $w=v$. 
\epv

Now:

\btm \label{cvx-3}
Let $h$ be a smooth metric on $S^2$. $h$ is induced on a convex surface
in $C^3_+$ if and only if $h$ is H-admissible. $h$ is induced on a
tamely convex surface if and only if $h$ is C-admissible. 
\etm

\bpv
Since $h$ is conformal to $\can_{S^2}$, so there exists a function
$u:S^2\rightarrow \R$ such that $e^{2u}\can_{S^2}$ is isometric to
$h$. If $P_0$ is any totally geodesic plane in $C^3+$, the graph $S$ of $u$
above $P_0$ has $h$ as its induced metric. Moreover, by the previous
lemma, $h$ is H-admissible if and only if $S$ has $B^*$ positive
definite, so if and only if $S$ is convex. And $h$ is C-admissible if
and only if $B^*$ has its eigenvalues in $(0,1)$, so if and only if $S$
is tamely convex. 
\epv

\begin{remark}
H-admissible metrics on $S^2$ have curvature $K<1$, while C-admissible
metrics on $S^2$ have curvature in $(-1,1)$. The converse,
however, is not true. 
\end{remark}

\bpv
Theorem \ref{cvx-3} shows that any H-convex metric is induced on a
convex surface in $C^3_+$, and lemma \ref{gauss} then indicates that it
has curvature strictly below $1$. Similarly C-convex metrics are induced
on tamely convex surfaces, which have curvature $K\in (-1,1)$ by lemma
\ref{gauss}.  
\epv

\section{Surfaces in $H^3$}

We will use in this section the results concerning the metrics on convex
hypersurfaces to understand the dual metrics on H-convex spheres in
$H^n$, and then on equivariant hypersurfaces.

\pg{Compact surfaces in $H^n, n\geq 4$}
As a consequence of theorems \ref{isom-n} and \ref{cvx-n}, we have for
$n\geq 4$: 

\btm \label{iin}
Let $h$ be smooth metric on $S^{n-1}$. $h$ is the horospherical metric
$I^*$ of a H-convex sphere $S$ in $H^n$ if and only if: 
\begin{itemize}
\item $h$ is conformal to $\can_{S^{n-1}}$;
\item $h$ is H-admissible, in the sense that it is conformal to
$\can_{S^{n-1}}$ and  
$ 2\ric_h - \frac{S_h}{n-2} - (n-3) h $ is everywhere negative
definite. 
\end{itemize}
$S$ is then unique up to the isometries of $H^n$. Moreover, $H$ is
tamely convex if and only if, at each point, 
all eigenvalues of $2(n-2)\ric_h-S_h h$ are in
$(-(n-2)(n-3), (n-2)(n-3))$.
\etm

\pg{Compact surfaces in $H^3$}
The same theorem holds in $H^3$ with the adequate notion of H-convexity;
it is a consequence of theorems \ref{isom-3} and \ref{cvx-3}. 

\btm \label{ii3}
Let $h$ be a smooth metric on $S^2$. It is the horospherical metric
$I^*$ of a H-convex immersed sphere $S$ in $H^3$ if and only if it is
H-admissible. It is the horospherical metric of a convex embedded sphere
$S\subset H^3$ if and only if it is C-admissible. In each case, $S$ is
unique up to the global isometries of $H^3$.  
\etm

\pg{Equivariant surfaces}
We consider now a surface $\Sigma$ of genus at least $2$. First we
introduce a class of metrics on $\Sigma$ in the following way -- this
definition is also a lemma. 

\bdf \label{df-hcvx-g}
Let $h$ be a smooth metric on $\Sigma$; we also call $h$ the pull-back
metric on the universal cover $\Sigmat$ of $\Sigma$. For each $x\in
\Sigmat$, there is a unique function $u_x:\Sigmat\rightarrow \R$ such that
$e^{-2u_x}h$ is isometric to a hemisphere of $(S^2, \can)$, and that
$u_x(x)=du_x(x)=0$. $h$ is {\bf H-admissible} if, for each $x$, the
hessian of $u_x$ is positive definite at $x$. 
$h$ is {\bf C-admissible} if, for each $x$, all eigenvalues of the
hessian at $x$ of $u_x$ are in $(0,1)$. 
\edf

Note that, here again, H-admissible metrics have $K<1$, and C-admissible
metrics have $K\in (-1,1)$, while the converse is false. 

\bpv
We only have to prove the existence and uniqueness of $u_x$. 

It is well known that there exists a unique hyperbolic metric in the
conformal class of $h$, i.e. a unique function $u:\Sigma\rightarrow \R$
such that $e^{-2u}h$ has constant curvature $-1$. We also call $u$ the
induced function on $\Sigmat$. Then $(\Sigmat, e^{-2u}h)$ is isometric to
$H^2$, and this defines a function $u$ on $H^2$ which is invariant under
an action of $\pi_1(\Sigma)$ by isometries. 

Now choose a dual plane $P_0\subset C^3_+$. Its induced metric is
isometric to that of $H^2$; choose an isometry between $P_0$ and 
$(\Sigmat, e^{-2u}h)$. This defines a function $u$ on $P_0$, and by
construction and lemma \ref{intrinsic}, the graph of $u$ above $P_0$ is
isometric to $(\Sigmat, h)$. We identify $\Sigmat$ with this graph. 

Now choose $x\in \Sigmat$, and let $P_1$ be the totally geodesic plane
tangent to $\Sigmat$ at $x$. $P_0$ is a graph above an hemisphere
$P_{1,+}$ of $P_1$, thus $\Sigmat$ is also the graph above $P_{1,+}$ of
a function $v$; by construction, $v$ satisfies the conditions on $u_x$.

Conversely, if $w$ is a function satisfying those conditions, then the
graph of $-w$ above $\Sigmat$ is a hemisphere of a totally geodesic plane
which is tangent to $\Sigmat$ at $x$, so $w=v$. 
\epv

This leads to a characterization of the metrics induced on equivariant
surfaces in $H^3$ as follows.

\btm \label{ie}
A smooth metric $h$ on $\Sigma$ is the horospherical metric of a
H-convex equivariant immersion whose representation fixes a plane if and
only if $h$ is H-admissible. It is the horospherical metric of a convex
equivariant embedding whose representation fixes a plane if and
only if $h$ is C-admissible. The equivariant immersion/embedding is then
unique up to global isometries. 
\etm

\bpv 
First note that any metric $h$ on $\Sigma$ has an equivariant isometric
embedding into $C^3_+$ whose representation fixes a dual plane. Indeed,
there is a unique function $u:\Sigma\rightarrow \R$ such that $e^{-2u}h$
is hyperbolic; $u$ can then be identified with an equivariant function
defined on a dual plane $P_0\subset C^3_+$, and then $(\Sigmat, h)$ is
isometric to the graph of $u$ aove $P_0$. 

The previous proof then indicates that $\Sigmat\subset C^3_+$ is convex
if and only if $h$ is H-admissible, and tamely convex if and only if $h$
is C-admissible. Therefore, the dual immersion in $H^3$ is H-convex if
and only if $h$ is H-admissible, and convex if and only if $h$ is
C-admissible. In this last case, the convexity implies that the
immersion is an embedding. 
\epv

Some kind of analoguous results in higher dimension might hold, but they
could be
less interesting since the metrics obtained are conformally flat, which
is a fairly strong condition. On the other hand they might be used to
put special (e.g. hyperbolic) metrics on conformally flat manifolds,
through deformations of equivariant sub-manifolds of $H^n$ or $C^n_+$.

\section{Hyperbolic manifolds with boundaries}

\pg{Why do all this ?}
As pointed out in the introduction, a natural question along conjecture
\ref{thurston} is to find the right boundary condition necessary to
obtain a unique hyperbolic metric on a given 3-manifold with
boundary. While conjecture \ref{thurston} strongly suggests that one should
consider the metric induced on the boundary, theorem \ref{tiie}
indicates that the third fundamental form of the boundary could be
another choice. 

The same question can be asked in higher dimensions, with hyperbolic
metrics replaced by Einstein metrics of negative curvature. A basic step
is taken in \cite{ecb}, where it is shown that any small deformation of
the metric induced on the boundary of a hyperbolic ball can be
``followed'' by an (essentially unique) Einstein deformation of the
metric in the interior. However, in this case again it is not completely
obvious whether the induced metric on the boundary is the right object
to consider. 

It appears clearly from recent work (see
e.g. \cite{graham-lee,graham-witten,witten,anderson-L2,anderson-close}) 
that, when one 
considers complete, conformally compact manifolds instead of metrics for
which the boundary is at finite distance, then the conformal class of
the boundary is what one needs. This does not indicate in any clear way
what one should use when the boundary is at finite distance, because, in
a conformally compact manifold, the hypersurfaces which are ``close'' to
the boundary in the conformal compact model are ``almost umbilical'', so
that the conformal class of the induced metric is also (asymptotically)
the conformal class of the second or third fundamental forms.  

So here again the solution advocated here is that the horospherical
metric might be the right thing to consider; the main argument is that,
for hyperbolic 3-manifolds, one can then obtain a satisfying existence
and uniqueness result in a very simple way. Of course the real challenge
will be to obtain similar results in higher dimensions, for Einstein
manifolds with boundary, or in other settings.

\pg{H-Admissible metrics}
We consider now a geometrically finite 3-manifold with boundary $(M, \dr
M)$ which admits a complete convex co-compact hyperbolic metric. Then
the universal cover 
of $(M, \dr M)$ is $(B^3, S^2)$, where $B^3$ is the 3-dimensional
ball. Moreover, if $h$ is a Riemannian metric on $\dr M$, then $h$
defines a complete metric on an open dense subset of $S^2$, which is invariant
under a conformal action of $\pi_1M$. Moreover, its conformal structure
defines a conformal structure on an open dense set of $S^2$, which
extends to a  
conformal metric on $S^2$, and the universal cover $\widetilde{\dr M}$ of
$\dr M$ has a unique conformal embedding into $S^2$ whose image is an
open dense set (see e.g. \cite{ahlfors}). 

We can thus define a proper class of metrics
on $\dr M$.

\bdf \label{h-adm-bord}
Let $h$ be a smooth metric on $\dr M$, and let $x\in \widetilde{\dr M}$. There
exists a unique function $u_x$ on $\widetilde{\dr M}$ such that $e^{-2u_x}
h$ extends to a constant curvature $1$ metric on $S^2$, and such that
$u_x(x)=du_x(x)=0$. We say that $h$ is {\bf H-admissible} if, for all
$x$, the hessian of $u_x$ at $x$ is positive definite. $h$ is {\bf
C-admissible} if, for all $x$, the eigenvalues of the  hessian of $u_x$
at $x$ are in $(0,1)$. 
\edf

Note that this definition coincides with definitions \ref{df-adm} and
\ref{df-hcvx-g} in the corresponding special cases. Again as above,
H-admissible metrics have curvature $K<1$, and C-admissible metrics
have curvature $K\in (-1,1)$. 

\bpv 
Again we only have to prove the existence and uniqueness of $u_x$. 

We know that there exists a function $u$ on the universal cover of $\dr
M$ such that $e^{-2u}h$ is isometric to an open dense subset of
$S^2$. This defines $(\Sigmat, h)$ as the graph of $u$ above an open
dense subset of a totally
geodesic plane $P_0$ (with the induced metric). 

The rest of the proof is just like for definition \ref{df-adm} and
\ref{df-hcvx-g}, and uses the uniqueness of the conformal change of metric. 
\epv

\pg{Existence and uniqueness}
We can now state the analog of conjecture \ref{thurston} for the
horospherical metric. 

\btm \label{thm-bord}
Let $h$ be a smooth metric on $\dr M$. 
\begin{enumerate}
\item $h$ is the horospherical metric
of a H-convex immersion $\phi$ of $\dr M$ in $M$ for a complete
hyperbolic metric  $g$ on $M$, such that the image of $\phi$ can be deformed
through immersions to the boundary at infinity of $M$, if and only if
$h$ is H-admissible. $g$ and $\phi$ are then unique.  
\item $h$ is the horospherical metric of $\dr M$ for a hyperbolic metric
$g$ on $M$, such that $\dr M$ is convex, if and only if $h$ is
C-admissible. $g$ is then unique.
\end{enumerate}
\etm

\bpv
We already know from the proof of definition \ref{h-adm-bord} that
$(\widetilde{\dr M}, h)$ is isometric to the graph of a unique (up to global
isometries) graph above a totally geodesic plane $P_0$. 

Taking the dual surface
in $H^3$ gives an immersion $\phi$ of $\widetilde{\dr M}$ in $H^3$ which
is H-convex if 
$h$ is H-admissible, and convex if $h$ is C-admissible. 

Moreover, $\pi_1M$ acts by
conformal transformations on $P_0$, so, by lemma \ref{isometries}, by
isometries on $H^3$. By construction, $\phi(\widetilde{\dr M})$ is
invariant under those isometries. Thus $(\dr M, h)$ is isometric
to the quotient by $\pi_1M$ of the image of $\phi$ with its horospheric
metric.  
\epv

Note that, if $h$ is only H-convex, we only obtain a priori an
immersion of $\dr M$ in $M$, which can be deformed through immersions to
an embedding. If $h$ is C-admissible, on the other hand, $\dr M$ is
obtained as a convex surface in $M$, so it is embedded (and it bounds a
convex domain in $M$). 

It should be pointed out that theorem \ref{ie} is a direct consequence
of theorem \ref{thm-bord}. Indeed consider the manifold $(M,\dr
M)=(\Sigma\times [-1,1], \Sigma\times \{ -1,1\})$, where $\Sigma$ is a
surface of genus at least $2$, and take on $\dr M$ a metric which is
identical on both copies of $\Sigma$. Then the uniqueness statement in
theorem \ref{thm-bord} implies that the metric $g$ obtained will have a
$\Z/2\Z$ symmetry, which will exchange the two connected components of
$\dr M$. Therefore $\widetilde{\dr M}$ will be immersed/embedded in $H^3$ as
two equivariant surfaces, symmetric with respect to a plane which is
fixed by both representations.

\pg{Higher dimensions}
Similar results might hold in the corresponding cases in higher
dimension, with conformally flat metrics on the boundary. This should
not be too interesting, however, since conformally flat metrics should
be quite rigid in this situation.

\section{Moreover}

\pg{An elementary approach}
A large part of what we have described here can be reduced essentially to a
simple (but remarkable) property. Let $H$ be a complete oriented
hypersurface in 
$H^n$, which is "uniformly H-convex" in the most natural sense. Let $u$
be a function on $h$, with a differential which is "small". For each point
$x\in H$, consider the horosphere $h_x$ tangent to $H$ at $x$, and its
equidistant horosphere $h_x'$ at distance $u(x)$. Then let $H'$ be the
envelope of the horospheres $h_x'$, and let $\phi$ be the map sending
$x\in H$ to the point $\phi(x)\in H'$ where $h_x'$ is tangent to $H'$
(this is well defined if $u$ and $H$ are well behaved. Then $\phi$ is an
isometry between $(H, e^{2u}I_H^*)$ and $(H', I_{H'}^*)$. 

Of course this is basically a translation, in purely $H^n$ terms, of the
basic properties of the metric on $C^n_+$, as described in lemma
\ref{intrinsic}. Moreover the statement is quite imprecise concerning
the precise conditions on $u$; of course things are clear in $C^n_+$,
the point is only that $u$ has to be such that the graph of $u$ above
$H^*$ (which will be the dual of $H'$) remains convex, so that $H'$
remains H-convex. More generally, I guess that some of the results
obtained here could be achieved without using $C^n_+$, but I doubt
whether it could improve the clarity of this matter.

\pg{Symmetric spaces and dualities}
Given a symmetric space $G/K$, there is a quite general way of
constructing other spaces (of the form $G/H$, for various choices of
$H\subset G$) which are in "duality" with $G/K$ -- see e.g. 
\cite{helgason-gga,helgason-gass}. The duality between $H^n$ and $C^n_+$
can be seen as a special case of this (with  $G=\SO(n,1)$, $K=\SO(n)$
and $H=\isom(\R^{n-1})$), just like the duality between
$H^n$ and $S^n_1$ (with $G=\SO(n,1)$, $K=\SO(n)$ and
$H=\SO(n-1,1)$). In this general setting there is a natural -- and well
understood -- duality between the functions or distributions on a space and
on its dual. The duality between the hypersurfaces can be put in this
context by replacing a hypersurface by some measure which it defines,
the dual hypersurface is then the support of the dual measure. 

In the
cases which we have described, however, one should not use the measure
associated to the area form on the hypersurfaces, since the duality
would then act with a factor equal to the Gauss-Kronecker curvature of
the hypersurfaces (in the case of the $H^n/S^n_1$ duality) or the
determinant of $E+B$ (in the $H^n/C^n_+$ duality). Rather one should
normalize this area measure by a factor $1/\sqrt{\det(B)}$ or
$1/\sqrt{\det(E+B)}$ in $H^n$, and $1/\sqrt{\det(B^*)}$ in $S^n_1$ or
$C^n_+$. 

A natural question is to understand to what extend the duality
properties of hypersurfaces in those spaces extend from the cases
described above to a more general setting, and what one could get out of
it.

\pg{Induced metrics and third fundamental form}
One striking feature of the results above is that they are simpler to
obtain -- and more powerful in some cases -- that the corresponding
results obtained for convex (hyper-)surfaces when one considers on them
the induced metrics or third fundamental forms. This leads to the idea
that those results could be used as a tool to obtain results on the
induced metrics or third fundamental form; for this one should obtain
rigidity results on the way the induced metric (resp. third fundamental
form) varies when a deformation changes the horospherical metric.

\pg{Einstein manifolds, etc}
The most natural framework in which conjecture \ref{thurston} could be
extended is the theory of negatively curved Einstein manifolds with
boundaries; indeed, in dimension 3, negatively curved Einstein metrics
are the same as hyperbolic metrics. 

An elementary (and far too restrictred) first step was taken in this
direction in \cite{ecb} (see also \cite{sem-era,sem} for some strikingly
related rigidity results obtained by very different methods). 
The outstanding problem there, however, is that the infinitesimal
rigidity result which is needed -- stating that an infinitesimal
deformation of the interior metric induces a non-trivial deformation of
the boundary metric -- is only obtained when the boundary is umbilical. 

A natural question is therefore whether an analog of the horospherical
metric (maybe defined as $I+2\II+\III$) could lead to some infinitesimal
rigidity result for Einstein manifolds with boundary; this would open
the door to possible results on the existence and/or uniqueness of
Einstein metrics inducing a given horospherical metric on the boundary. 

Note that the theory concerning complete metrics is rather more
advanced; in that case one only prescribes the conformal structure on
the boundary at infinity, and the Einstein metrics are required to be
conformally compact. In dimension 3 it is just the classical
Ahlfors-Bers theorem, while in higher dimension the theory seems to be
advancing (see the previous section for references).

\bibliographystyle{alpha}

\begin{thebibliography}{Lab92b}

\bibitem[Ahl66]{ahlfors}
L.~V. Ahlfors.
\newblock {\em Lectures on quasiconformal mappings}.
\newblock D. Van Nostrand Co., Inc., Toronto, Ont.-New York-London, 1966.
\newblock Manuscript prepared with the assistance of Clifford J. Earle, Jr. Van
  Nostrand Mathematical Studies, No. 10.

\bibitem[Ale51]{alex}
A.~D. Aleksandrov.
\newblock {\em Convex Polyhedra}.
\newblock GITTL, 1951.
\newblock Russian; english translation to appear, Springer.

\bibitem[Ale58]{Al}
A.~D. Aleksandrov.
\newblock {\em Vestnik Leningrad Univ.}, 13(1), 1958.

\bibitem[Anda]{anderson-close}
M.~T. Anderson.
\newblock Closedness of the space of {AHE} metrics on 4-manifolds.
\newblock math.DG/0012167.

\bibitem[Andb]{anderson-L2}
M.~T. Anderson.
\newblock {$L^2$} curvature and volume renormalization of {AHE} metrics on
  4-manifolds.
\newblock math.DG/0011051.

\bibitem[AZ67]{AZ}
A.~D. Aleksandrov and V.~A. Zalgaller.
\newblock {\em Intrinsic Geometry of Surfaces}, volume~15 of {\em Translations
  of Mathematical Monographs}.
\newblock AMS, 1967.

\bibitem[Bes87]{Be}
A.~Besse.
\newblock {\em Einstein Manifolds}.
\newblock Springer, 1987.

\bibitem[Car16]{cartan}
E.~Cartan.
\newblock La d{\'e}formation des hypersurfaces dans l'espace euclidien r{\'e}el {\`a}
  n dimensions.
\newblock {\em Bull. Soc. Math. France}, 44:65--99, 1916.

\bibitem[GHL87]{GHL}
S.~Gallot, D.~Hulin, and J.~Lafontaine.
\newblock {\em Riemannian Geometry}.
\newblock Springer, 1987.

\bibitem[GL91]{graham-lee}
C.~R. Graham and J.~M. Lee.
\newblock Einstein metrics with prescribed conformal infinity on the ball.
\newblock {\em Adv. Math.}, 87:186--225, 1991.

\bibitem[Gro86]{PDR}
M.~Gromov.
\newblock {\em Partial Differential Relations}.
\newblock Springer, 1986.

\bibitem[GW]{graham-witten}
C.~R. Graham and E.~Witten.
\newblock Conformal anomaly of submanifold observables in {AdS/CFT}
  correspondence.
\newblock hep-th/9901021.

\bibitem[Hel94]{helgason-gass}
S.~Helgason.
\newblock {\em Geometric analysis on symmetric spaces}.
\newblock American Mathematical Society, Providence, RI, 1994.

\bibitem[Hel00]{helgason-gga}
S.~Helgason.
\newblock {\em Groups and geometric analysis}.
\newblock American Mathematical Society, Providence, RI, 2000.
\newblock Integral geometry, invariant differential operators, and spherical
  functions, Corrected reprint of the 1984 original.

\bibitem[HR93]{hr}
C.~D. Hodgson and I.~Rivin.
\newblock A characterization of compact convex polyhedra in hyperbolic 3-space.
\newblock {\em Invent. Math.}, 111:77--111, 1993.

\bibitem[Lab89]{L1}
F.~Labourie.
\newblock Immersions isom{\'e}triques elliptiques et courbes
  pseudo-holo\-morphes.
\newblock {\em J. Differential Geom.}, 30:395--424, 1989.

\bibitem[Lab92a]{L4}
F.~Labourie.
\newblock M{\'e}triques prescrites sur le bord des vari{\'e}t{\'e}s hyperboliques de
  dimension 3.
\newblock {\em J. Differential Geom.}, 35:609--626, 1992.

\bibitem[Lab92b]{L5}
F.~Labourie.
\newblock Surfaces convexes dans l'espace hyperbolique et {CP1}-structures.
\newblock {\em J. London Math. Soc., II. Ser.}, 45:549--565, 1992.

\bibitem[LS00]{iie}
F.~Labourie and J.-M. Schlenker.
\newblock Surfaces convexes fuchsiennes dans les espaces lorentziens {\`a}
  courbure constante.
\newblock {\em Math. Annalen}, 316:465--483, 2000.

\bibitem[Nir53]{N}
L.~Nirenberg.
\newblock The {Weyl} and {Minkowski} problem in differential geometry in the
  large.
\newblock {\em Comm. Pure Appl. Math}, 6:337--394, 1953.

\bibitem[Pog73]{Po}
A.~V. Pogorelov.
\newblock {\em Extrinsic Geometry of Convex Surfaces}.
\newblock American Mathematical Society, 1973.
\newblock Translations of Mathematical Monographs. Vol. 35.

\bibitem[RH93]{rh}
I.~Rivin and C.~D. Hodgson.
\newblock A characterization of compact convex polyhedra in hyperbolic 3-space.
\newblock {\em Invent. Math.}, 111:77--111, 1993.

\bibitem[RS99]{sem-era}
I.~Rivin and J.-M. Schlenker.
\newblock The {Schl{\"a}fli} formula in {Einstein} manifolds with boundary.
\newblock {\em Electronic Research Announcements of the A.M.S.}, 5:18--23,
  1999.

\bibitem[RS00]{sem}
I.~Rivin and J.-M. Schlenker.
\newblock The {Schl{\"a}fli} formula and {Einstein} manifolds.
\newblock Preprint math.DG/0001176, 2000.

\bibitem[Sch96]{these}
J.-M. Schlenker.
\newblock Surfaces convexes dans des espaces lorentziens {\`a} courbure
  constante.
\newblock {\em Commun. Anal. and Geom.}, 4:285--331, 1996.

\bibitem[Sch98]{shu}
J.-M. Schlenker.
\newblock M{\'e}triques sur les poly{\`e}dres hyperboliques con\-vexes.
\newblock {\em J. Differential Geom.}, 48(2):323--405, 1998.

\bibitem[Sch00]{dap}
J.-M. Schlenker.
\newblock Dihedral angles of convex polyhedra.
\newblock {\em Discrete Comput. Geom.}, 23(3):409--417, 2000.

\bibitem[Sch01]{ecb}
J.-M. Schlenker.
\newblock Einstein manifolds with convex boundaries.
\newblock Pr{\'e}\-publi\-ca\-tion no 98-12, Universit{\'e} de Paris-Sud;
  Commentarii Mathematici Helvetici, to appear, 2001.

\bibitem[Wit98]{witten}
E.~Witten.
\newblock Anti de {Sitter} space and holography.
\newblock {\em Adv. Theor. Math. Phys.}, 2:253--291, 1998.

\end{thebibliography}

\end{document}